\documentclass[twoside]{article}
\usepackage[a4paper]{geometry}
\usepackage[latin1]{inputenc} 
\usepackage[T1]{fontenc} 
\usepackage{RR}
\usepackage{hyperref}

\usepackage{amsmath,epsfig,graphics,graphicx}
\usepackage{amsfonts,amssymb}
\usepackage{color}
\usepackage{mathrsfs}
\usepackage{enumitem}
\usepackage{array}
\usepackage{float}
\usepackage{verbatim}

\definecolor{darkgreen}{RGB}{0,142,128}
\definecolor{darkgray}{RGB}{102,102,102}

\newcommand{\vpar}{v_\parallel}

\newcommand{\Bstar}{B_{\|s}^{\ast}}
\newcommand{\bstar}{\vec{b}_s^{\ast}}
\newcommand{\vecbstar}{\vec{b}_s^{\ast}}

\newcommand{\vecx}{\vec{x}}
\newcommand{\vecb}{\vec{b}}

\newcommand{\Bparstar}{B_{||s}^*}

\newcommand{\gradxi}{\vec{\nabla} x^i}

\newcommand{\gradB}{\vec{\nabla} B}
\newcommand{\gradF}{\vec{\nabla} F}
\newcommand{\gradG}{\vec{\nabla} G}

\newcommand{\gyroPhi}{\bar{\phi}}
\newcommand{\gradgyroPhi}{\vec{\nabla} \gyroPhi}

\newcommand{\Bgradxi}{\vec{B}\cdot\gradxi}

\newcommand{\muJgradxi}{\mu_0\vec{J}\cdot\gradxi}

\newcommand{\vgcEgradxi}{\vec{v}_{Es_{GC}}\cdot\gradxi}

\newcommand{\vDgradxi}{\vec{v}_{Ds}\cdot\gradxi}

\newcommand{\vgcEgradB}{\vec{v}_{Es_{GC}}\cdot\gradB}

\newcommand{\rhostar}{\rho^{\ast}}

\newcommand{\der}{\delta{}r}
\newcommand{\deth}{\delta{}\theta}

\let\pa\partial
\def\be{\begin{equation}}
\def\ee{\end{equation}}
\def\rad{r}
\def\brho{\tilde{\rho}}

\RRNo{8054}
\RRdate{September 2012}
\RRauthor{
G. Latu\thanks[sfn]{CEA Cadarache, IRFM bat. 513, 13108 Saint-Paul-les-Durance}%
  \and
V. Grandgirard\thanksref{sfn}%
\and J. Abiteboul\thanksref{sfn} 
\and M. Bergot\thanks[ca]{University of Strasbourg, 7 rue Descartes, 67000 Strasbourg  \& INRIA CALVI}
 \and N. Crouseilles\thanks{INRIA IPSO \& IRMAR, 263 av Général Leclerc, 35042 Rennes}
\and X. Garbet\thanksref{sfn} \and P. Ghendrih\thanksref{sfn} \and M. Mehrenberger\thanksref{ca}
\and Y. Sarazin\thanksref{sfn} \and H. Sellama\thanksref{ca} \and E. Sonnendrücker\thanksref{ca} 
\and D. Zarzoso\thanksref{sfn}
}
\authorhead{Latu \& al}
\RRtitle{Précision du mouvement de particules non perturbées dans un code gyrocinétique}
\RRetitle{Accuracy of unperturbed motion of particles in a gyrokinetic semi-Lagrangian code}
\titlehead{Unperturbed motion of particles}
\RRresume{Une description imparfaite de l'équilibre peut conduire à
  des artéfacts numériques dans des simulations turbulentes
  gyrocinétiques. Aussi, le solveur de Vlasov and le schéma
  d'intégration en temps ont un impact fort sur la conservation de
  quantités physiques, notamment lorsque les simulations sont en temps
  long. L'équilibre et le solveur Vlasov doivent être finement choisis
  et paramétrés pour conserver les états constants (équilibre) et pour
  autoriser de bonnes propriétés de conservation en temps (en
  commençant par la conservation de la masse). Plusieurs cas tests
  illustratifs sont donnés pour montrer les problèmes numériques
  typiques que l'on peut observer si les choix pris ne sont pas
  adéquats. Nous expliquons pourquoi le schéma FSL (Forward
  Semi-Lagrangian) apporte une réponse à certains problèmes. Des
  simulations en configuration cylindrique et torique de GYSELA sont
  présentées qui utilisent FSL.}

\RRabstract{Inaccurate description of the equilibrium can yield to
  spurious effects in gyrokinetic turbulence simulations. Also, the
  Vlasov solver and time integration schemes impact the conservation
  of physical quantities, especially in long-term
  simulations. Equilibrium and Vlasov solver have to be tuned in order
  to preserve constant states (equilibrium) and to provide good
  conservation property along time (mass to begin with). Several
  illustrative simple test cases are given to show typical spurious
  effects that one can observes for poor settings. We explain why
  Forward Semi-Lagrangian scheme bring us some benefits. Some toroidal
  and cylindrical GYSELA runs are shown that use FSL.}
\RRmotcle{gyrocinétique, lois de conservation, FSL}
\RRkeyword{gyrokinetic, conservation laws, FSL}
\RRprojets{CALVI \& IPSO}
\RCNancy 

\begin{document}
\makeRR   

\section{Introduction}
Inaccurate description of the gyrokinetic equilibrium can yield 
unphysical excitation of zonal flow
oscillations\cite{Ang06}. Moreover, as stated in \cite{Ido08}, it is
important to define the initial condition using a relevant gyrokinetic
equilibrium, especially in collisionless full-f simulations. In the
following, accuracy aspects are investigated for both the gyrokinetic
equilibrium and the Vlasov solver used in the GYSELA code. If proper
care is not taken for Vlasov solver and gyrokinetic initial
equilibrium, one can observe that some conservation properties are not
satisfied, for example mass conservation.

The gyrokinetic framework of the study is explained in the next
section. Then, quite experimental investigations are shown to
understand numerics associated with equilibrium and mass
conservation. Forward Semi-Lagrangian (FSL) scheme is presented
briefly. The splitting in a linear and a non-linear part of the Vlasov
equation is explained.  Results of numerical simulation for 4D
cylindrical test cases and toroidal test cases are presented in the
last section. Forward Semi-Lagrangian scheme helps a lot to preserve
good accuracy on mass and total energy.

\section{Description of the context}
\subsection{Gyrokinetic Vlasov equation}
The time evolution of the gyro-center coordinates $(\vecx,\vpar,\mu)$
of species $s$ is given by the collision-less electrostatic
gyrokinetic equations:
\begin{align}
  \frac{{\rm d}x^i}{\rm dt} &= \vpar \vecbstar\cdot\gradxi + \vgcEgradxi + \vDgradxi \label{dxidt}\\
  m_s\frac{{\rm d}\vpar}{\rm dt} &= -\mu \vecbstar\cdot\gradB - e_s \vecbstar\cdot\gradgyroPhi + \frac{m_s\vpar}{B}\vgcEgradB \label{dvpardt}
\end{align}
where $x^i$ corresponds to the $i$-th covariant coordinate of $\vecx$,
B is the magnetic field and $\Bstar$ and $\vecbstar$ are defined as:
\begin{align}
  \Bstar &= B+\frac{m_s\vpar}{e_sB}\mu_0\vecb\cdot\vec{J} \label{Bstar} \\
  \bstar &= \frac{\vec{B}}{\Bstar}+\frac{m_s\vpar}{e_s\Bstar}\,\frac{\mu_0\vec{J}}{B} \label{vecbstar}
\end{align}
The advection terms are:
\begin{align}
  \vecbstar\cdot\gradxi & = b_s^{\ast i} =\frac{\Bgradxi}{\Bstar}+\frac{m_s\vpar}{e_s\Bstar} \frac{\muJgradxi}{B} \label{vpar_drift}\\
  \vDgradxi &= v_{Ds}^i = \left(\frac{m_s \vpar^2+\mu B}{e_s\Bstar\,B}\right)[B,x^i]\label{vD_drift}\\
  \vgcEgradxi &= v_{Es_{GC}}^i =\frac{1}{\Bstar}[\gyroPhi,x^i]\label{vE_drift}\\
  \vgcEgradB  &= -\frac{1}{\Bstar}[B,\gyroPhi]\label{vgcEgradB_drift}
\end{align}
the Poisson bracket is defined by
$\scriptstyle[F,G]\ =\ \vecb\cdot\left(\gradF\times
\gradG\right)$. The velocity parallel to the magnetic field is
$\vpar$. The magnetic moment \mbox{$\mu=m_s v_\perp^2/(2B)$} is an
adiabatic invariant with $v_\perp$ the velocity in the plane
orthogonal to the magnetic field. Vacuum permittivity is denoted
$\mu_0$.  The term $v_{E_{GC}}$ represents the electric $E\times B$
drift velocity of the guiding-centers and $v_D$ the curvature drift
velocity.
The Jacobian in phase space is $J_s(r,\theta)\Bparstar(r,\theta,\vpar)$.

Some references concerning the framework we use to solve this equation can
be found in \cite{virginie,vlasovia06,braeunig11,sarazin11,abiteboul11}.
\subsection{Quasineutrality equation}

The quasineutrality equation and parallel Amp\`{e}re's law close the
self-consistent gyrokinetic Vlasov-Maxwell system. However, in an
electrostatic code, the field solver reduces to the numerical solving
of a Poisson-like equation\,\cite{hahm}.  In tokamak configurations,
the plasma quasineutrality (denoted QN) approximation is currently
assumed\,\cite{virginie}.  Electron inertia is ignored, which means
that an adiabatic response of electrons is supposed. We define the
operator
$\textstyle{\nabla\!}_{\perp}=(\pa_r,\frac{1}{r}\pa_{\theta})$. We
note $n_0$ the equilibrium density, $B_0$ the magnetic field at the
magnetic axis and $T_e(r)$ the electronic temperature. We have also
$B(r,\theta)$ the magnetic field, $J_0$ the Bessel function of first
order and $k_{\perp}$ the transverse component of the wave vector.
Hence, the QN equation can be written in dimensionless
variables\\ \scalebox{.80}{
\begin{minipage}[b]{1.25\textwidth}\centering
\be
-\frac{1}{n_0(r)}\nabla_{\perp}\,.\left[\frac{n_0(r)}{B_0}\nabla_{\perp}\Phi(\rad,\theta,\varphi)\right]+\frac{1}{T_e(r)}\left[\Phi(r,\theta,\varphi)-\left<\Phi\right>_{FS}(r)\right]
= \brho(r,\theta,\varphi)
\label{qneq}
\ee
\end{minipage}
}

\medskip
\noindent where $\brho$ is defined by\\
\scalebox{.80}{
\begin{minipage}[b]{1.25\textwidth}\centering
\be
\label{brho}
\brho(r,\theta,\varphi)=\frac{2\pi}{n_0(r)}\int \Bparstar d\mu\int
dv_{\parallel}
J_0(k_{\perp}\sqrt{2\mu})(\bar{f}-\bar{f}_{eq})(r,\theta,\varphi,v_{\parallel},\mu).
\ee 
\end{minipage}
}
with $\bar{f}_{eq}$ representing local ion Maxwellian equilibrium,
and $\left<.\right>_{FS}(r)$ the average on the flux surface labelled by $r$.

\section{Unperturbed motion of particles}
\subsection{Simplification of the equations}
Let suppose that a collisionless equilibrium distribution function
$f_{eq}$ is taken for the initial distribution function:
$f_{t=0}(r,\theta,\varphi,v_\parallel)$. By definition, if the Vlasov
and Quasi-neutrality equations are accurately solved, the distribution
function $f_t$ should remain constant over time, \textit{i.e.}  equal
to $f_{eq}$. Also, the electric potential $\Phi(r,\theta,\varphi)$
should remain equal to zero because the right hand side of the
Quasi-neutrality equation, the integral of $f-f_{eq}$, is zero.\\ Then,
one quantitative measure of Vlasov solver accuracy is its capability
to preserve an equilibrium distribution function over time, assuming
$\Phi=0$. In order to get a simple tool to evaluate the Vlasov solver,
one can also restrict to the case $\mu=0$.  It is worthwhile rewriting
gyrokinetic equations in this reduced setting:
\begin{align}
  \frac{{\rm d}x^i}{\rm dt} &= \vpar \vecbstar\cdot\gradxi + \vDgradxi \label{dxidt2}\\
  m_s\frac{{\rm d}\vpar}{\rm dt} &= 0 \label{dvpardt2}
\end{align}
It gives
\begin{align}
  \frac{{\rm d}x^i}{\rm dt} &= \frac{\Bgradxi}{\Bstar}\vpar{}+\frac{m_s\vpar{}^2}{e_s\Bstar} \frac{\muJgradxi}{B}+\left(\frac{m_s \vpar^2}{e_s\Bstar\,B}\right)[B,x^i]\\
  \frac{{\rm d}\vpar}{\rm dt} &= 0
\end{align}
and finally, assuming that the current J is perpendicular to the poloidal plane,
\begin{align}
  \frac{{\rm d}r}{\rm dt} &= \frac{B_r}{\Bstar}\vpar{}+\frac{m_s
    \vpar^2}{e_s\Bstar\,B}[B,r] \label{drdt} \\
  \frac{{\rm d}\theta}{\rm dt} &=
  \frac{B_\theta}{\Bstar}\vpar{}+\frac{m_s
    \vpar^2}{e_s\Bstar\,B}[B,\theta] \label{dthetadt} \\
  \frac{{\rm d}\varphi}{\rm dt} &=
  \frac{B_\varphi}{\Bstar}\vpar{}+\frac{m_s\vpar{}^2}{e_s\Bstar}
  \frac{J_\varphi}{B}+\frac{m_s
    \vpar^2}{e_s\Bstar\,B}[B,\varphi] \label{dphidt} 
\end{align}
Particles keep a constant parallel velocity over time.  Then the
trajectory of one particle in phase space remains in the 3D space
dimensions.

\subsection{Trajectory of particles}
\label{vlasovsplit}
The GYSELA code uses a Semi-Lagrangian scheme combined with a Strang
splitting. To guarantee second order accuracy in time, the
following sequence\,\cite{virginie} is used:
$(\hat{\vpar}/2,\hat{\varphi}/2,\hat{r\theta},\hat{\varphi}/2,\hat{\vpar}/2)$,
where factor 1/2 denotes a shift over $\Delta{}t/2$. In this sequence,
$\textstyle\hat{\vpar}$ and $\textstyle\hat{\varphi}$ shifts
correspond to 1D advections along $\vpar$ and $\varphi$ directions,
and $\textstyle\hat{r\theta}$ represents a 2D advection (displacement
in $r$ and $\theta$ are strongly coupled). Let us now focus on the
discretization of the 2D advection which requires specific attention.

\paragraph{Taylor expansion}
For a single $\textstyle\hat{r\theta}$ advection in the
semi-Lagrangian approach, a displacement $\textstyle(\der,\deth)$ is
computed to find the foot of the characteristics that ends at a grid
node $(r_i,\theta_j)$. This displacement is computed thanks to a
Taylor expansion of the fields acting on
particles\,\cite{virginie}. One can derive two fields $\alpha$ and
$\beta$ from Eqs. $(\ref{drdt},\ref{dthetadt},\ref{dphidt})$ such as:
$$\begin{pmatrix}\der{}\\\deth{}\end{pmatrix}_{(r_i,\theta_j)}=\alpha{}(r_i,\theta_j)\Delta{}t+\beta{}(r_i,\theta_j)\Delta{}t^2+O(\Delta{}t^3)$$
 The problem is that, if $\der{}$ or $\deth{}$ is larger than grid element
 sizes denoted $\Delta{}r$ and $\Delta{}\theta$, then this approximation can be
 really inaccurate.

\paragraph{Precomputed trajectories}
Another method is described here to find out the foot of the
characteristics associated to each grid node. Since the fields acting
on particles, Eqs. (\ref{dxidt2},\ref{dvpardt2}), do not depends on
time $t$, one can approximate
$\textstyle(\der{},\deth{})_{(r_i,\theta_j)}$ once for all. We have used
Runge-Kutta time integration scheme (RK2) with a small time step to
precompute these trajectories. This approach is possible only for
the linear terms of Eqs. (\ref{dxidt2},\ref{dvpardt2}), but not for
non-linear terms, present in Eqs. (\ref{dxidt},\ref{dvpardt}), that
depend on $E_t$ and $\Phi_t$. Let us choose a $\delta{}t$, such as
$M\,\delta{}t = \Delta{}t$ with $M\in\mathbb{N}$ and $M$ large enough. One can build a series $r^n$
and $\theta^n$ as follows:
$$\begin{pmatrix}r^{n+\frac{1}{2}}\\\theta^{n+\frac{1}{2}}\end{pmatrix}=
\begin{pmatrix}r^{n}\\\theta^{n}\end{pmatrix}+\frac{\delta{}t}{2}\,\alpha{}(r^n,\theta^n)$$
$$\begin{pmatrix}r^{n+1}\\\theta^{n+1}\end{pmatrix}=
\begin{pmatrix}r^{n}\\\theta^{n}\end{pmatrix}+\delta{}t\ \alpha{}(r^{n+\frac{1}{2}},\theta^{n+\frac{1}{2}})$$
The initial condition is set to $r^0=r_i$, $\theta^0=\theta_j$. After $M$ steps of
Runge-Kutta iteration, it gives 
$$\begin{pmatrix}\der{}\\\deth{}\end{pmatrix}_{(r_i,\theta_j)}=\begin{pmatrix}r^{M}\\\theta^{M}\end{pmatrix}$$
Because these trajectories can be computed only once when the simulation
starts, $M$ can be taken quite large (we will assume $M=64$ in the
following). These precomputations will not impact significantly the
global simulation time. Other time integration schemes have been
tried: RK3, RK4, and also larger values of $M$, but these
modifications does not impact significantly the precision in our test
cases.

\subsubsection{Forward Semi-Lagrangian scheme}

The Forward Semi-Lagrangian scheme\,\cite{fsl09} (the acronym is FSL)
is an alternative to the classical Backward Semi-Lagrangian scheme
(BSL). It is based on forward integration of the characteristics. The
distribution function is updated on an eulerian grid, and the
pseudo-particles located on the mesh nodes follow the
characteristics of the equation forward for one time step, and are
then deposited on nearest nodes. While the main cost of the BSL method comes
from the interpolation step, FSL spends most of computational time in
the deposition step. The FSL scheme can be set up such as mass is
preserved along time.

We present the approach taken for the 2D advection in
$(r,\theta)$. In our setting, the algorithm of Forward Semi-Lagrangian is the
following (with $S$ the cubic B-spline):
\begin{itemize}
\item Step 1: Compute the spline coefficients $w_{k,l}^n$ ($f^n$ is
  known) such that
$$  (J_s \Bparstar f^{n})_{i,j} = \sum_{k,l} w^n_{k,l} S(x_i-x_k) S(y_j-y_l)$$
\item Step 2: Integrate forward in time the characteristics from $t^n$ to
  $t^{n+1}$, given as initial data the grid points $(x_k,y_l)$. We
  obtain the new particle positions  $(x^\star_k,y^\star_l)$ at time
  $t^{n+1}$ by following the characteristics.
\item Step 3: Deposition of particles localized at $(x^\star_k,y^\star_l)$
$$ f^{n+1}_{i,j} = \left(\sum_{k,l} w^n_{k,l} S(x_i-x^\star_k) S(y_j-y^\star_l)\right) / (J_s \Bparstar)_{i,j} $$
\end{itemize} 
\subsection{Numerical experiments}

\subsubsection{Accuracy of conservation of invariants}

An equilibrium solution of the collisionless gyrokinetic equation must
satisfy some conditions. To get an equilibrium for the Vlasov
equation, it suffices to take an arbitrary function of constants of
motion in the unperturbed characteristics. In an axisymmetric toroidal
configuration, a gyrokinetic Vlasov equilibrium is defined by three
constants of motion: the magnetic moment~$\mu$, the energy
\mbox{$\scriptstyle \mathcal{E}=\vpar^2/2+\mu\,B(r,\theta)$}, and the canonical
toroidal angular momentum\,\cite{abiteboul11} \mbox{$\scriptstyle
  P_\varphi=\psi(r)+I\vpar{}/B(r,\theta)$}. The $\psi(r)$ function is defined
thanks to the safety factor $q(r)$ by the relation: \mbox{$\scriptstyle
  d\psi/dr=-B_0\,r/q(r)$}.

To get a good accuracy in integrating the Vlasov equation in time, and then to 
conserve exactly an initial equilibrium distribution function, we have
to take care of the correctness of:
\begin{itemize}
\item the initial equilibrium setting; distribution function will be a
  function of the constants of motion: $\mu, \mathcal{E}, P_\varphi$;
\item computing the displacements in the Vlasov solver (finding the
  foot of the characteristic); a sufficiently small time step has to
  be chosen for the time integration scheme for example, but also the
  resolution in space must be sufficient to have enough accuracy to
  interpolate displacement fields;
\item interpolate the distribution function in the Vlasov solver with
  enough resolution; the number of points in the computational domain
  has to be tuned for this issue;
\end{itemize}
We will be able to estimate quantitatively the quality of the
distribution function after several time steps by:
\begin{itemize}
\item measuring the difference between the initial equilibrium
  function and the distribution function at a given time step.
\end{itemize}

%
%


\subsubsection{First test case: equilibrium conservation}
\label{firstcase}
The initial distribution function is set up to
\begin{equation}
f_{t=0}(r,\theta,\varphi,v_\parallel)=
\mathbf{1}_{[P_{\varphi{}1},P_{\varphi{}2}]}
(P_\varphi-P_{\varphi{}1})^m
(P_\varphi-P_{\varphi{}2})^m\,sin(\gamma{}\,P_\varphi)
\label{initf1}
\end{equation}
\begin{figure}[h]
\begin{minipage}[t]{.48\textwidth}
\includegraphics[trim=2mm 15mm 63mm 20mm,clip,width=\textwidth]{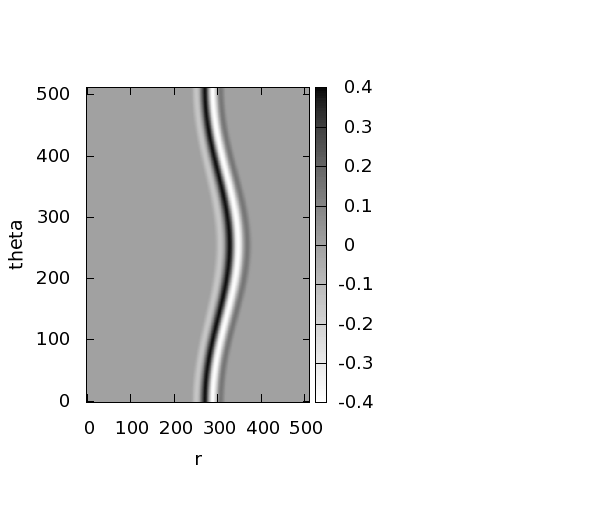}
\end{minipage}
\begin{minipage}[t]{.48\textwidth}
\includegraphics[trim=2mm 15mm 63mm 20mm,clip,width=\textwidth]{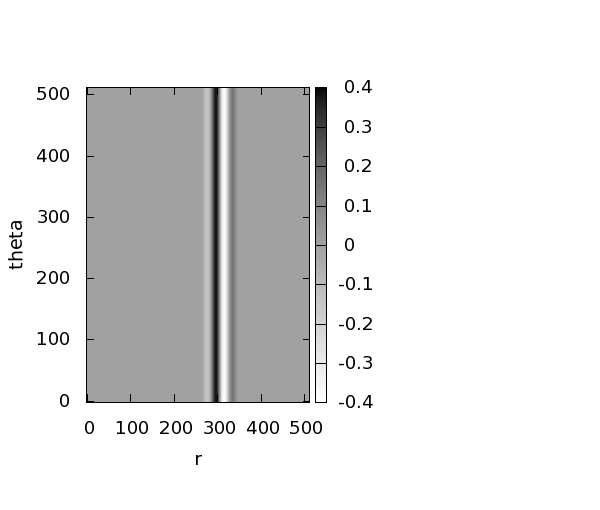}
\end{minipage}
\caption{Initial function $f_{t=0}$ at $v_{\parallel}\!=\!-3\,v_{th0},\varphi\!=\!0$ (left), $v_{\parallel}\!=\!0,\varphi\!=\!0$ (right)}
\label{distribinit}
\end{figure}
where the invariant $P_\varphi$ is a function of $r,\theta,\vpar$ and
$\mathbf{1}$ is the indicator function. We will look at the norms
\mbox{$u_p(t)=\|f_{t}-f_{t=0}\|_p$} with \mbox{$p=1,\infty$}. They are defined as
\begin{eqnarray*}
  u_\infty(t)=sup \{|f_{t}-f_{t=0}|:\forall (r,\theta,\varphi,\vpar)\}\\
  u_1(t)=\int |f_{t}-f_{t=0}| J_s \Bparstar dr d\theta d\varphi d\vpar
\end{eqnarray*}

 Let
us take the following parameters :
\mbox{$\gamma=4\pi/(P_{\varphi{}2}-P_{\varphi{}1})$, $m=1$,
  $\rhostar=0.01$}. The values $(P_{\varphi{}1}, P_{\varphi{}2})$ are
choosen adequately in order to have an interval $[P_{\varphi{}1},
  P_{\varphi{}2}] $ that is strictly included in the computational
domain
\mbox{$r\in[r_{min},r_{max}],\theta\in[0,2\pi],v_{\parallel}=\!-3\,v_{th0}$}. 
A poloidal cut of this initial function is sketched for two values of
parallel velocity in Figure \ref{distribinit}. As the initial
function depends only on the invariant $P_\varphi$, it should be
conserved by time integration of the Vlasov equation. The duration of
the experiment is denoted $t_{max}$. The timing unit is the ion
cyclotronic time $\Omega^{-1}$, it will also be the case in all plots
of this document. The Vlasov equation is solved during all time
steps with $\Phi$ forced to zero. An ideal simulation should preserve all norms
\mbox{$\forall p,\ \ u_p(t)=0$}. Practically and in the worst case,
the function $u_p(t)$ will move quickly away from 0 as $t$ grows.
%
%
\begin{figure}[H]
\begin{minipage}{.48\textwidth}
\includegraphics[trim=25mm 22mm 80mm 55mm,clip,width=\textwidth]{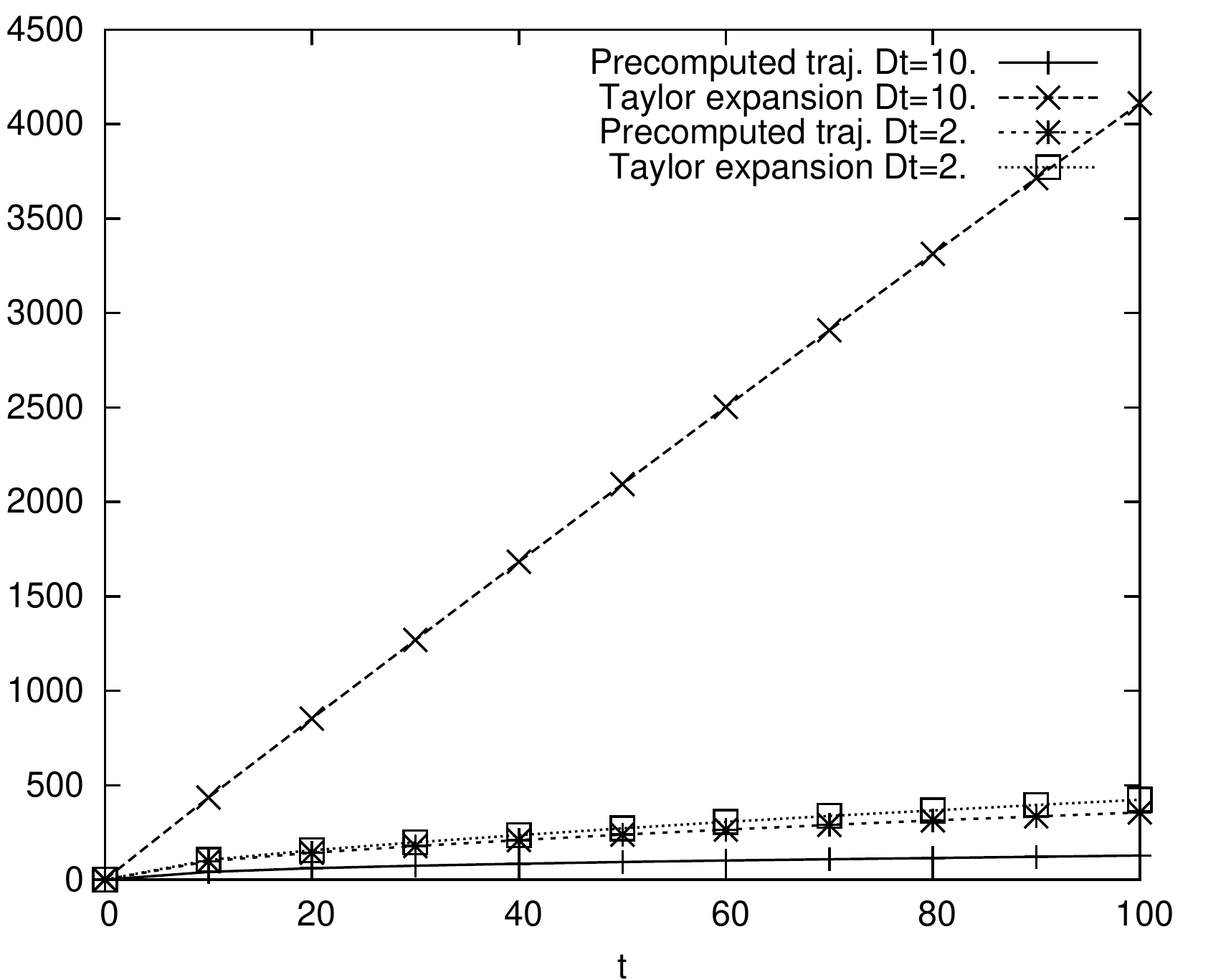}
\caption{Evolution of the $L_1$ norm $u_1(t)$ with $N_r=512,N_\theta=512$}
\label{scanDt1}
\end{minipage}~
\begin{minipage}{.48\textwidth}
\includegraphics[trim=25mm 22mm 80mm 55mm,clip,width=\textwidth]{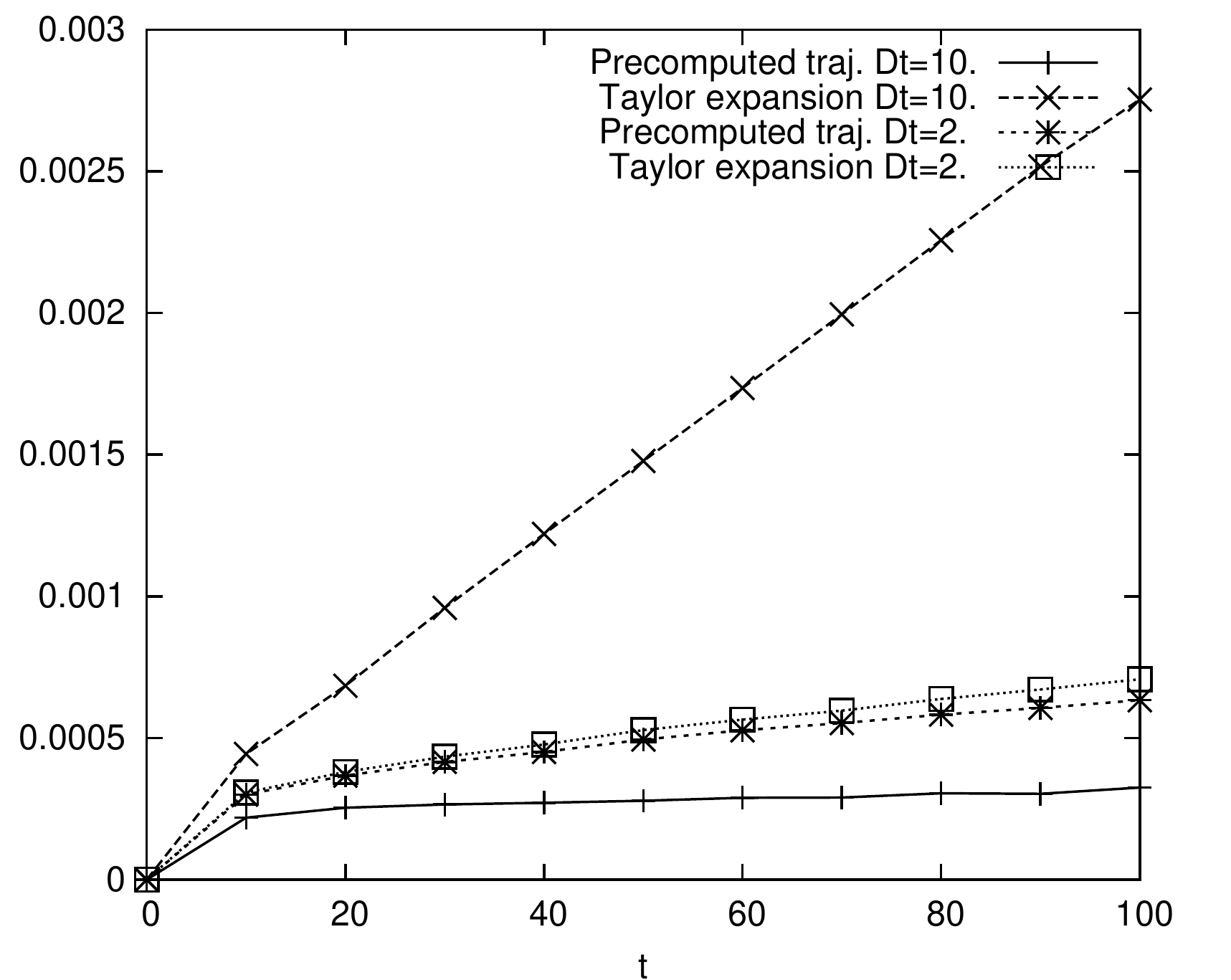}
\caption{Evolution of the $L_\infty$ norm $u_{\infty}(t)$ with $N_r=512,N_\theta=512$}
\label{scanDt3}
\end{minipage}
\end{figure}
\begin{figure}[H]
\begin{minipage}[t]{.48\textwidth}
\includegraphics[trim=2mm 15mm 63mm 20mm,clip,width=\textwidth]{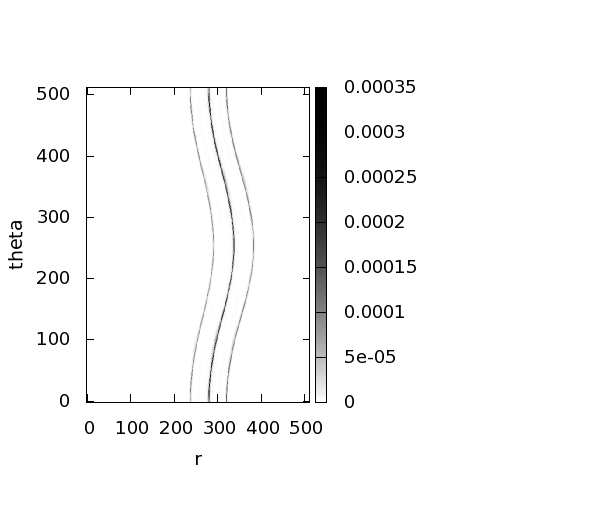}
\end{minipage}
\begin{minipage}[t]{.48\textwidth}
\includegraphics[trim=2mm 15mm 63mm 20mm,clip,width=\textwidth]{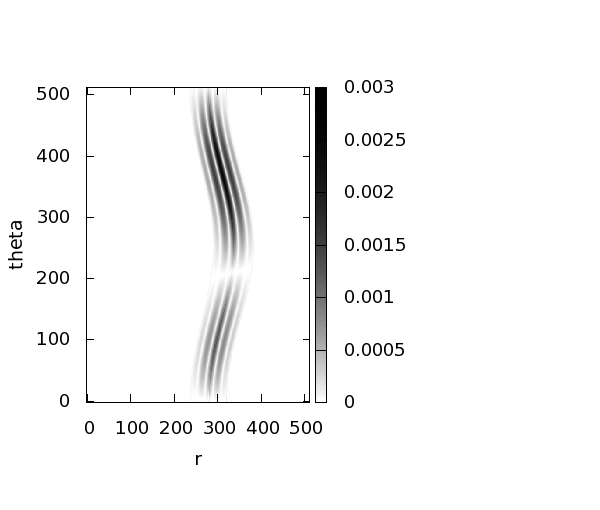}
\end{minipage}
\caption{distance between $f_t$ and $f_{t=0}$ (modulus of the difference) at $t\!=\!100$ with $\Delta{}t\!=\!10$ and \textit{Precomputed trajectories} (left), distance $\|f_t-f_{t=0}\|$ at $t\!=\!100$ with $\Delta{}t\!=\!10$ and \textit{Taylor expansion} (right)}
\label{distrib0}
\end{figure}

In the experiment of Figures \ref{scanDt1} and \ref{scanDt3}, the
$L_1, L_\infty$ norms are drawn for some simulations.  A comparison is
made between the two methods presented previously that are used to
calculate the feet of the characteristics. The parameter $\Delta{}t$
is set to $2$ and $10$, which is near typical values used in GYSELA
simulations and \mbox{$t_{max}\!=\!100$}. The computational domain is
quite refined in $(r,\theta)$ with $N_r\!=\!512,
N_\theta{}\!=\!512$. One can remark that \textit{Precomputed
  trajectories} strategy is more accurate than \textit{Taylor
  expansion}. Nevertheless, the difference shown here is significant
for $\Delta{}t=10$, but quite negligible for small $\Delta{}t=2$.  On
Figure \ref{distrib0}, the differences to the initial distribution
function at final state, are shown for an arbitrary value of
$v_{\parallel}=-3\,v_{th0}$. One can see here where errors are located
using both strategies: \textit{Precomputed trajectories} at the
center, and \textit{Taylor expansion} at right. One can see that the
amplitude of the error is larger in the case of \textit{Taylor
  Expansion}, and also less spatially localized.

\begin{figure}[H]
\begin{minipage}{.48\textwidth}
\includegraphics[trim=25mm 22mm 80mm 55mm,clip,width=\textwidth]{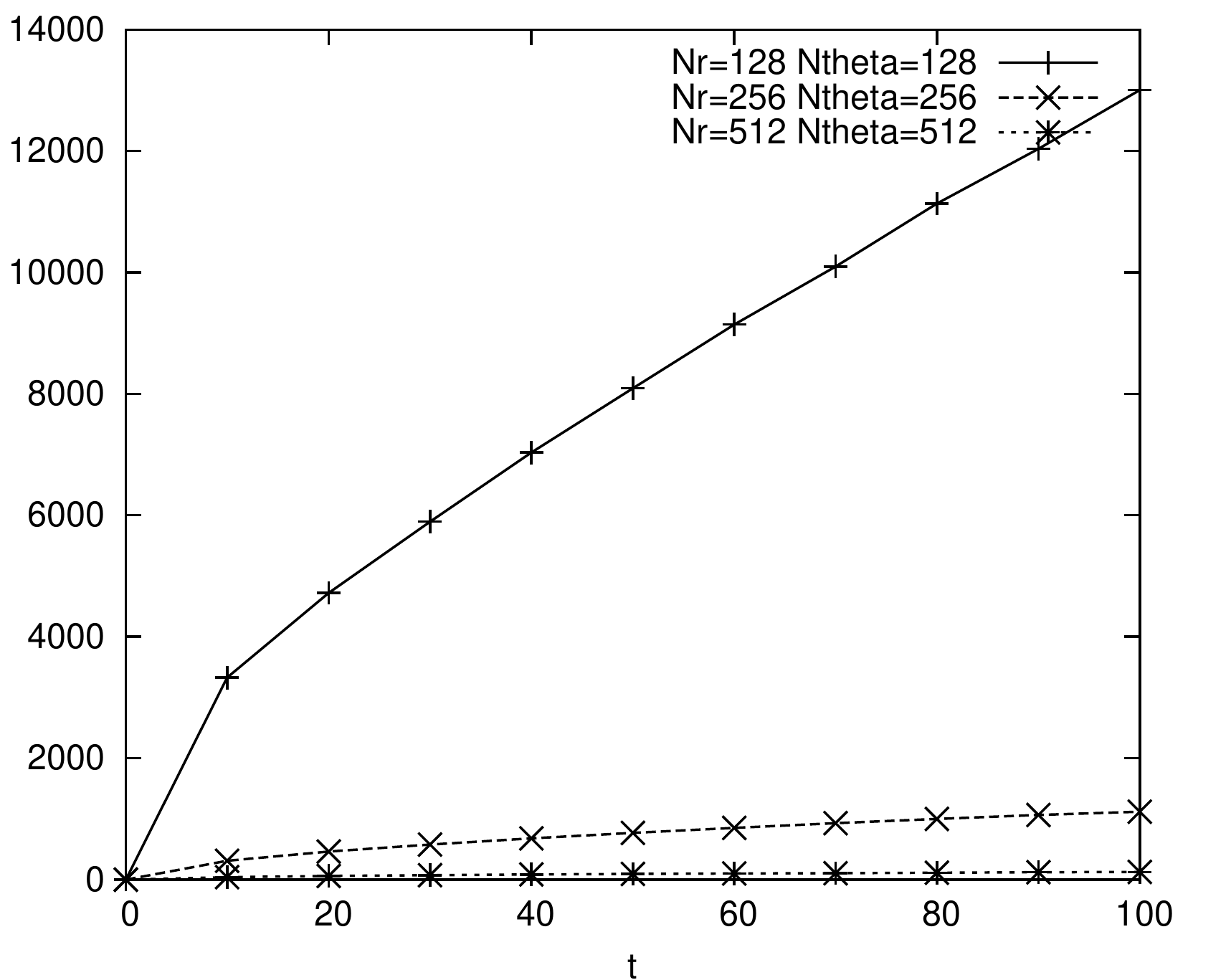}
\caption{Evolution of the $L_1$ norm $u_1(t)$  for different poloidal mesh sizes}
\label{scanSp1}
\end{minipage}~
\begin{minipage}{.48\textwidth}
\includegraphics[trim=25mm 22mm 80mm 55mm,clip,width=\textwidth]{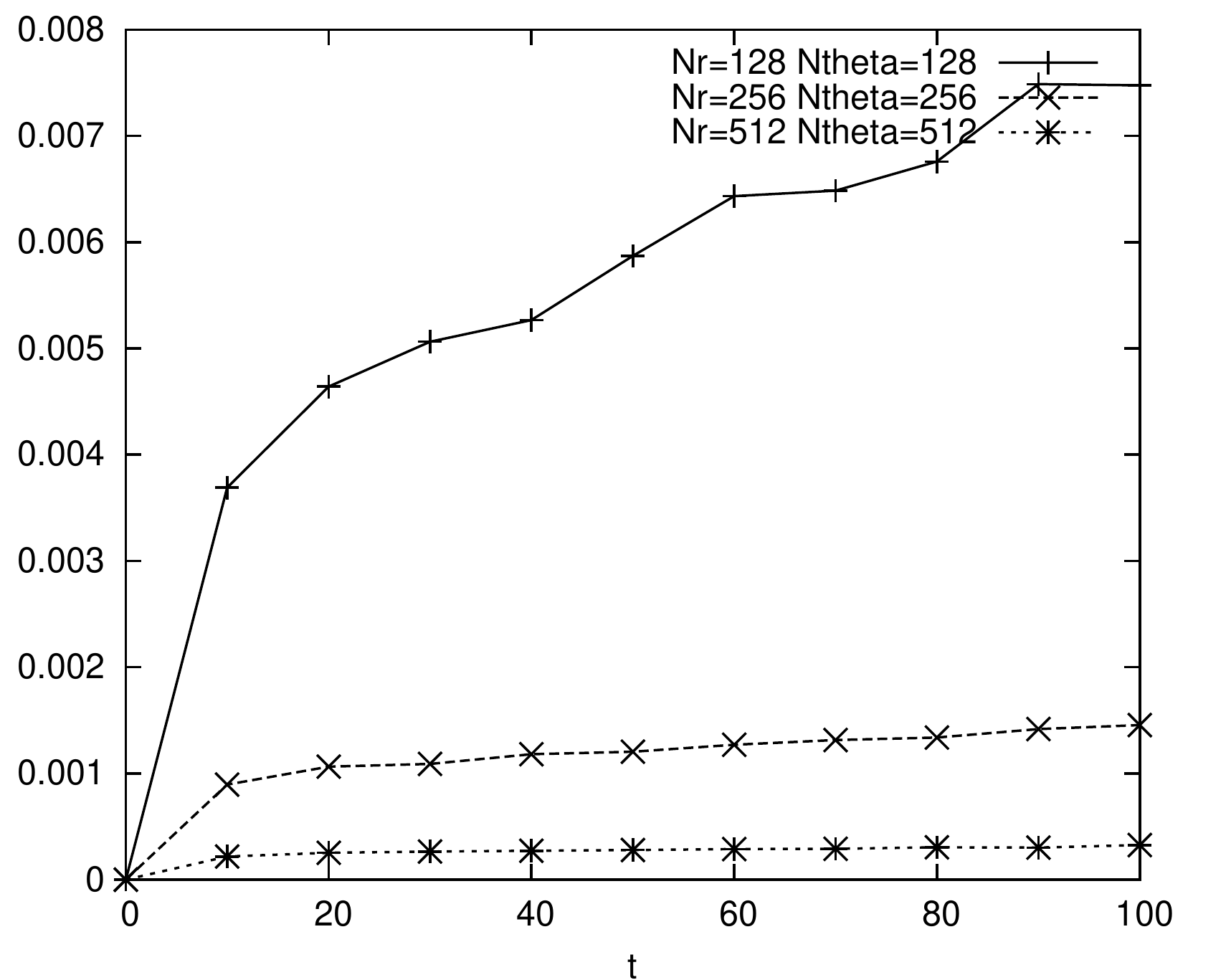}
\caption{Evolution of the $L_\infty$ norm $u_{\infty}(t)$ for different poloidal mesh sizes}
\label{scanSp3}
\end{minipage}
\end{figure}

On Figures \ref{scanSp1} and \ref{scanSp3}, some $L$-p norms produced
by three simulations are shown (\textit{Precomputed trajectories}
method is used). Between two consecutive simulations, the sizes of the domain
$(N_r,N_{\theta})$ have been multiplied by a factor two.  While
increasing $(N_r,N_{\theta})$, one can notice that simulations become
more and more accurate, \textit{i.e.} L-$p$ norms are reduced. The
increase rate of the error depends on the factor
$N_r\,N_\theta$.

\begin{figure}[h]
\begin{minipage}{.48\textwidth}
\includegraphics[trim=25mm 22mm 80mm 55mm,clip,width=\textwidth]{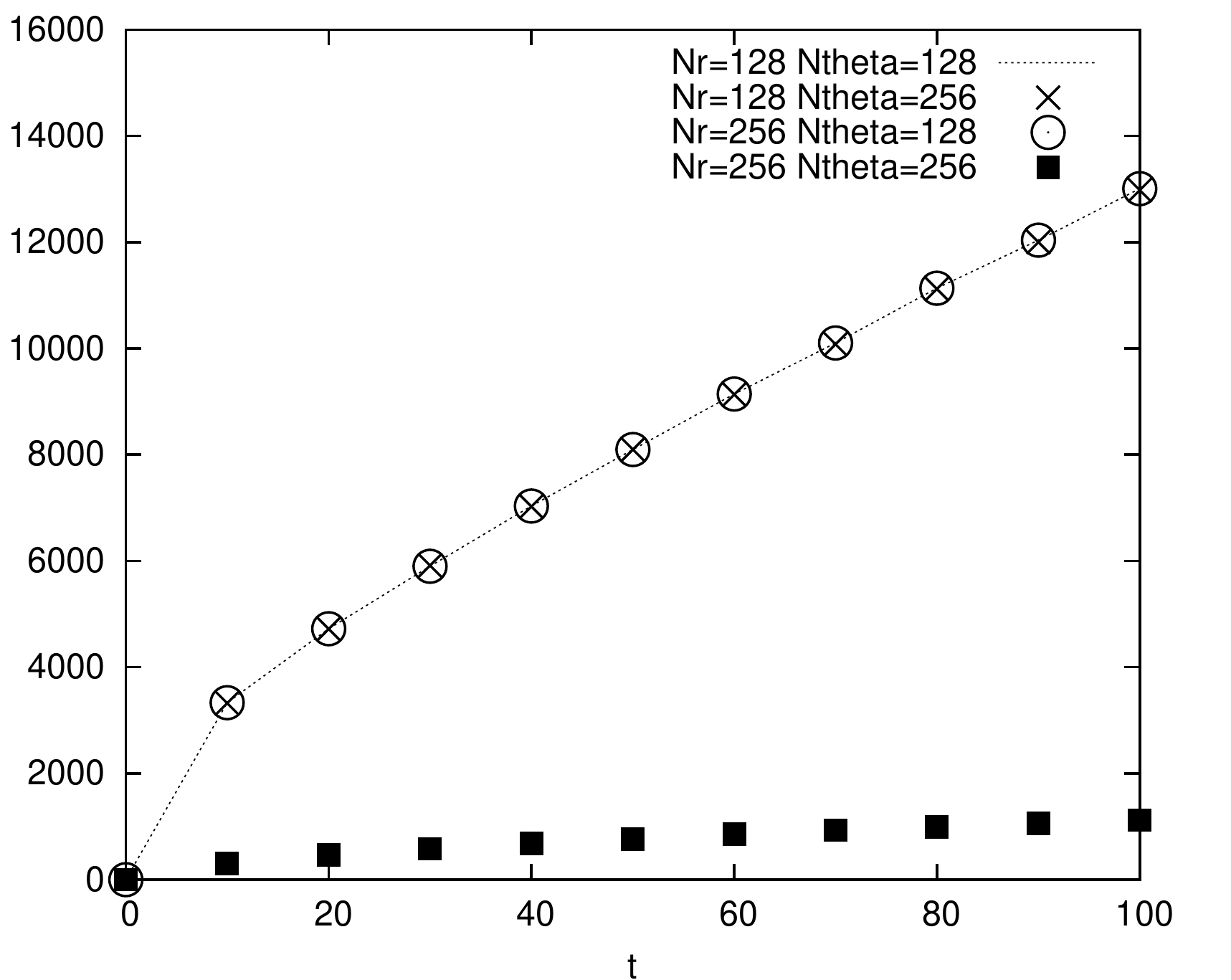}
\caption{Evolution of the $L_1$ norm $u_1(t)$  for different poloidal mesh sizes}
\label{scanSp4}
\end{minipage}~
\begin{minipage}{.48\textwidth}
\includegraphics[trim=25mm 22mm 80mm 55mm,clip,width=\textwidth]{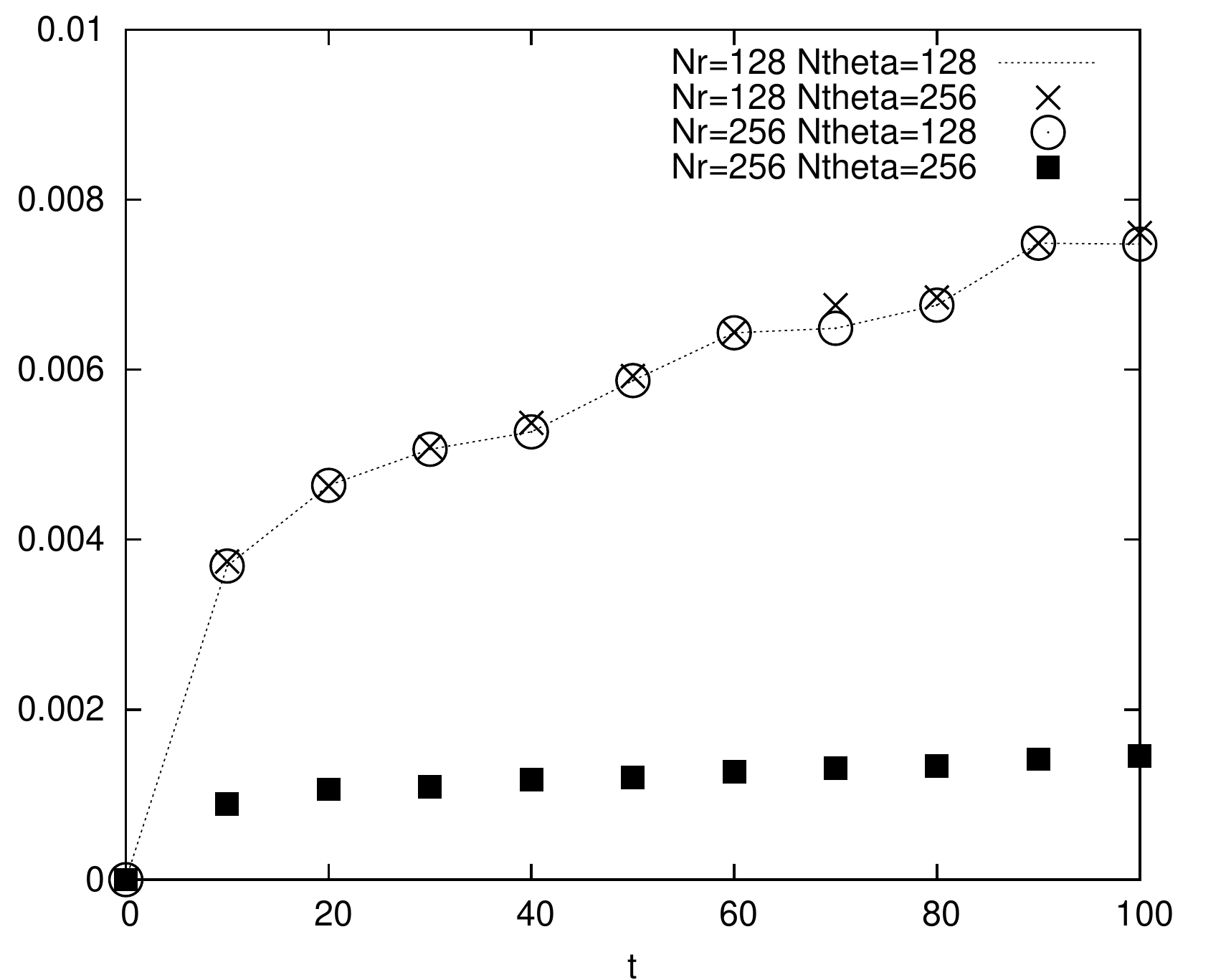}
\caption{Evolution of the $L_\infty$ norm $u_{\infty}(t)$ for different poloidal mesh sizes}
\label{scanSp6}
\end{minipage}
\end{figure}

The domain sizes $N_r$ and $N_{\theta}$ are increased independantly
each other in the Figures \ref{scanSp4} and \ref{scanSp6}
(\textit{Precomputed trajectories} method is used). It is clear that
the error is only of the order of the smallest size among $N_r$ and
$N_\theta$ sizes. As far as Equations (\ref{dxidt2},\ref{dvpardt2})
are concerned, it is then noteworthy that $N_r$ and $N_\theta$ should
be chosen close to each other in order to improve accuracy.

%

\subsubsection{Second test case: shear flows}
\label{secondcase}

Compared to the previous experiments, the initial distribution
function considered in the present section is taken with an
inhomogeneity in $\theta$. Only the \textit{Precomputed trajectories}
method has been used for these simulations, because it is more
accurate in any case. As variable $\theta$ is not an invariant, the
distribution function is expected to evolve in time. For the sake of
simplicity and to reduce computational costs, the simulation is setup
with $N_{\varphi}=1, N_{v_\parallel}=1$. The initial distribution
function is now
$$f_{t=0}(r,\theta,\varphi,v_\parallel)=
\mathbf{1}_{[P_{\varphi{}1},P_{\varphi{}2}]}
(P_\varphi-P_{\varphi{}1})^m (P_\varphi-P_{\varphi{}2})^m
\mathbf{1}_{[\theta{}_1,\theta{}_2]} (\theta-\theta{}_1)^m
(\theta-\theta{}_2)^m$$

\begin{figure}[H]
\begin{minipage}[t]{.32\textwidth}
\includegraphics[trim=0mm 15mm 94mm 30mm,clip,width=\textwidth]{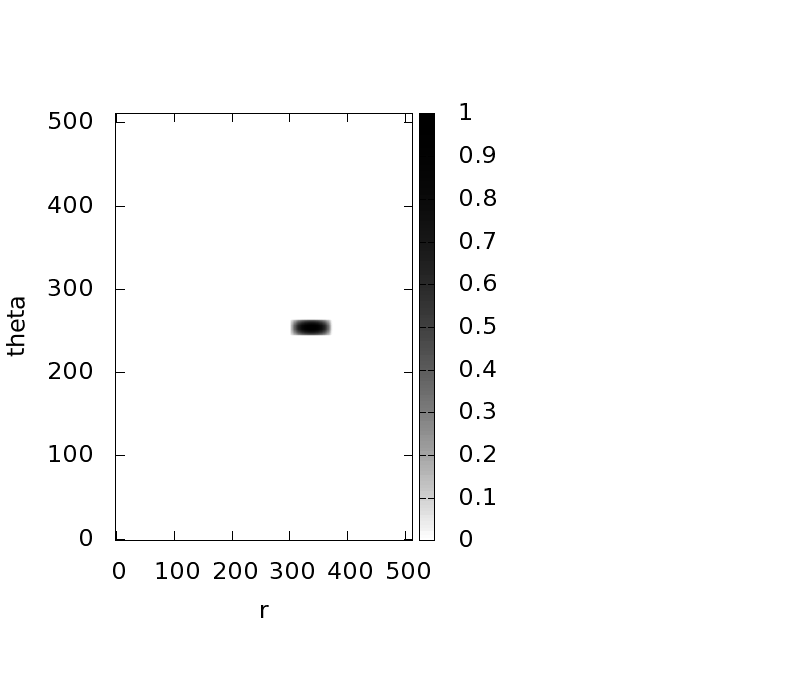}
\end{minipage}
\begin{minipage}[t]{.32\textwidth}
\includegraphics[trim=0mm 15mm 94mm 30mm,clip,width=\textwidth]{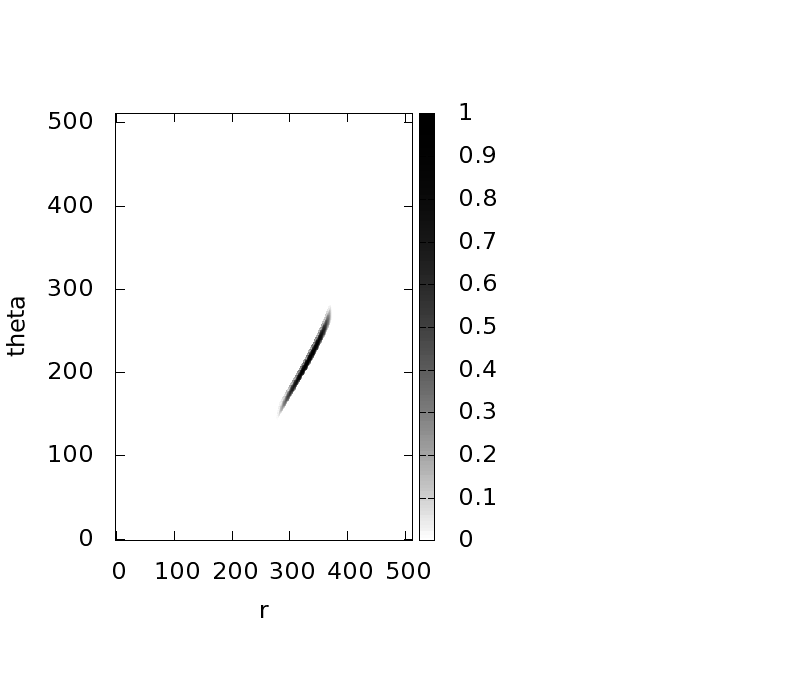}
\end{minipage}
\begin{minipage}[t]{.32\textwidth}
\includegraphics[trim=0mm 15mm 94mm 30mm,clip,width=\textwidth]{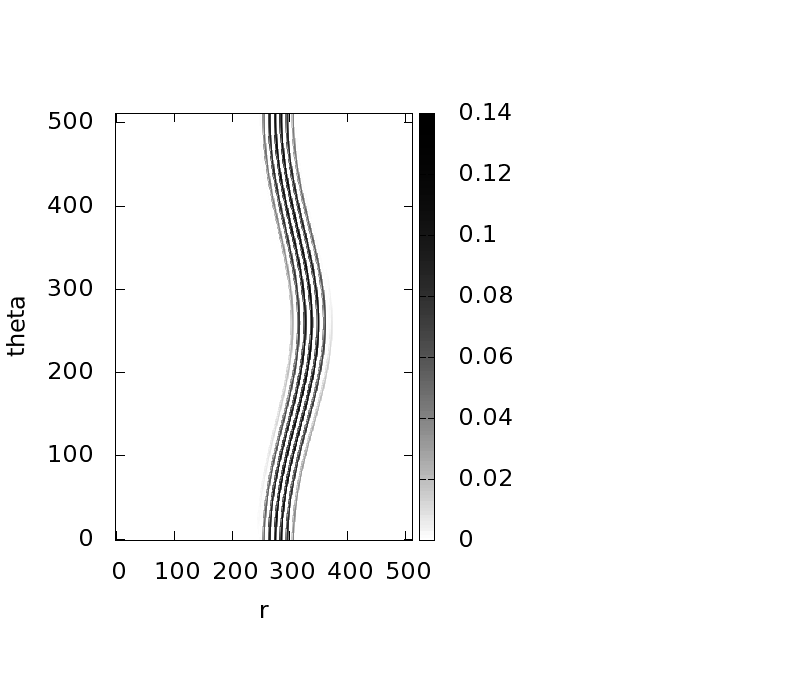} 
\end{minipage}
\caption{Initial function $f_{t=0}$ with $v_{\parallel}\!=\!-3\,v_{th0}, N_r=512, N_\theta=512$ (left), distribution function $f$ at \mbox{$t=2000$} with BSL scheme (middle), distribution function $f$ at \mbox{$t=50000$} with BSL scheme (right)}
\label{shear20}

\begin{minipage}[t]{.32\textwidth}
\includegraphics[trim=0mm 15mm 94mm 30mm,clip,width=\textwidth]{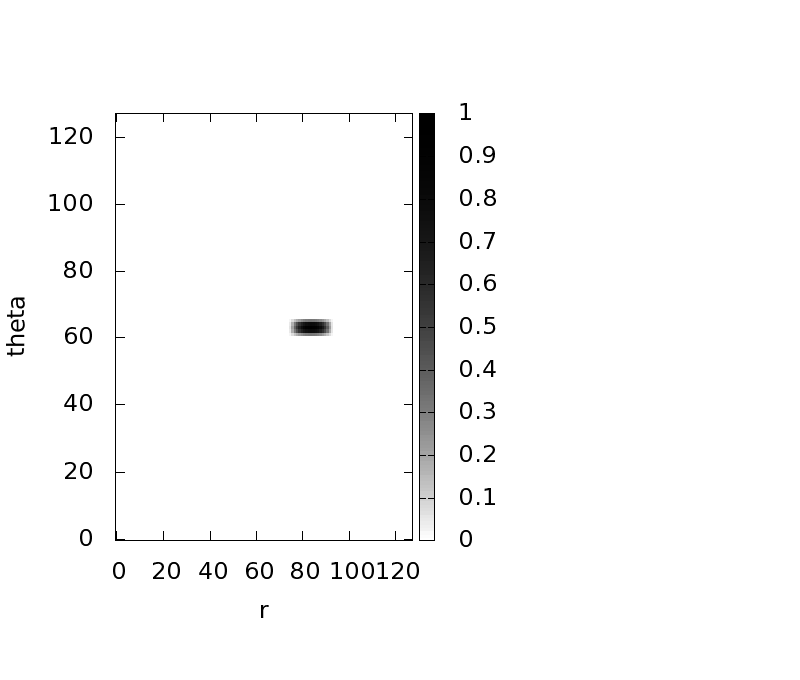}
\end{minipage}
\begin{minipage}[t]{.32\textwidth}
\includegraphics[trim=0mm 15mm 94mm 30mm,clip,width=\textwidth]{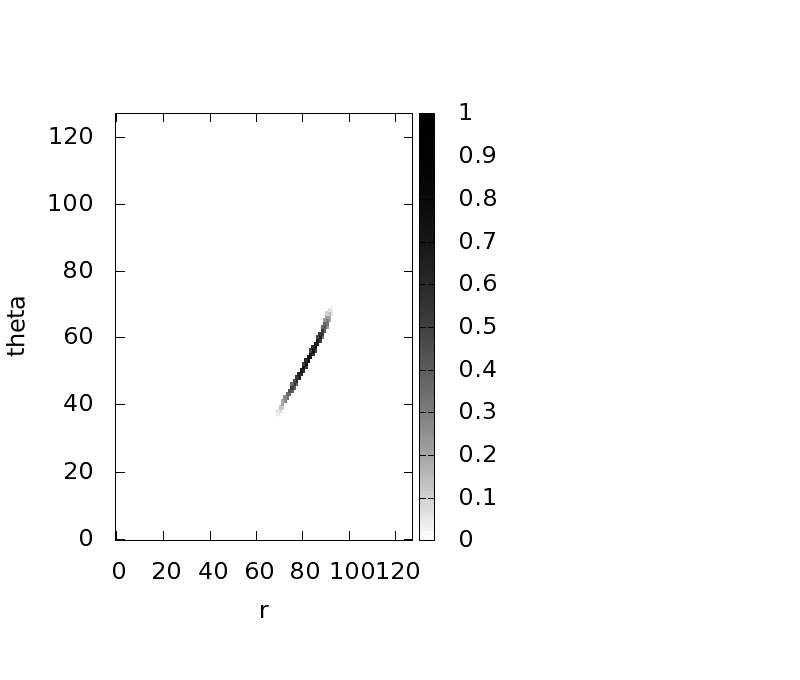}
\end{minipage}
\begin{minipage}[t]{.32\textwidth}
\includegraphics[trim=0mm 15mm 94mm 30mm,clip,width=\textwidth]{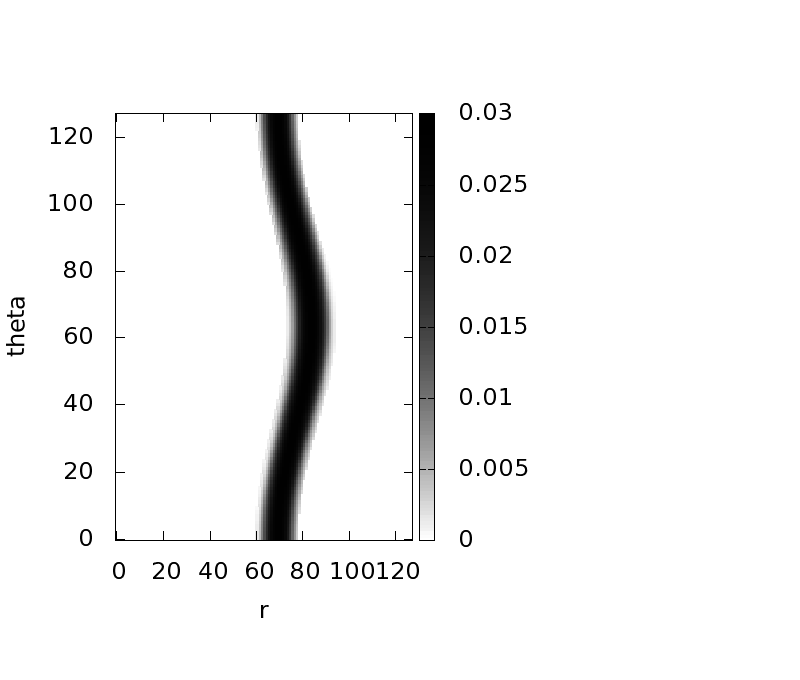}
\end{minipage}
\caption{Initial function $f_{t=0}$ with $v_{\parallel}\!=\!-3\,v_{th0}, N_r=128, N_\theta=128$ (left), distribution function $f$ at \mbox{$t=2000$} with BSL scheme (middle), distribution function $f$ with BSL scheme at \mbox{$t=50000$}  (right)}
\label{shear0}

\begin{minipage}[t]{.32\textwidth}
\includegraphics[trim=0mm 15mm 94mm 30mm,clip,width=\textwidth]{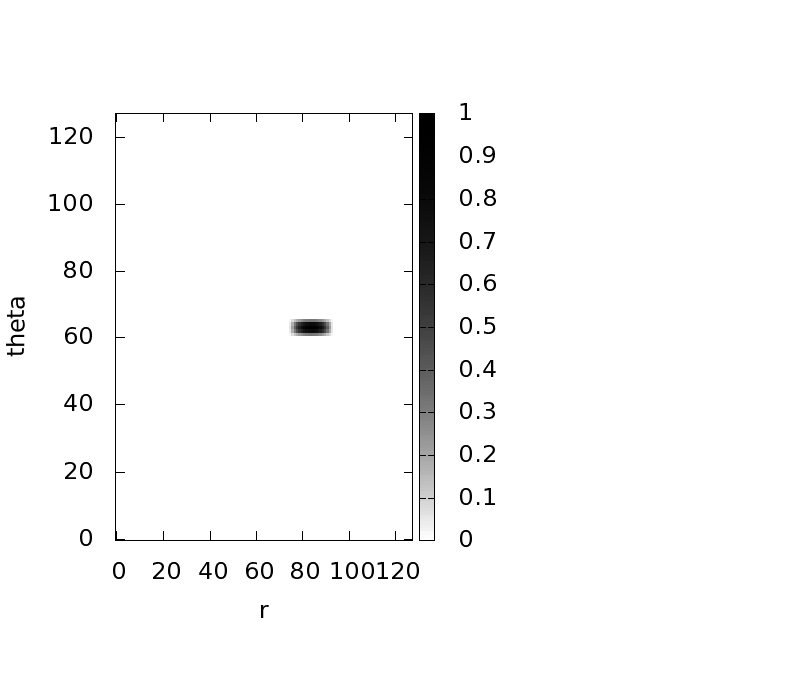}
\end{minipage}
\begin{minipage}[t]{.32\textwidth}
\includegraphics[trim=0mm 15mm 94mm 30mm,clip,width=\textwidth]{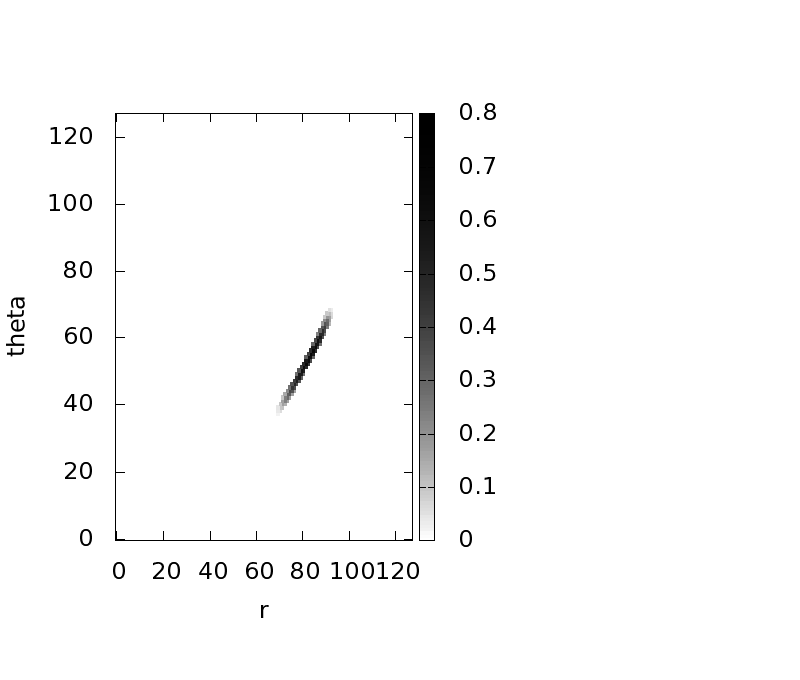}
\end{minipage}
\begin{minipage}[t]{.32\textwidth}
\includegraphics[trim=0mm 15mm 94mm 30mm,clip,width=\textwidth]{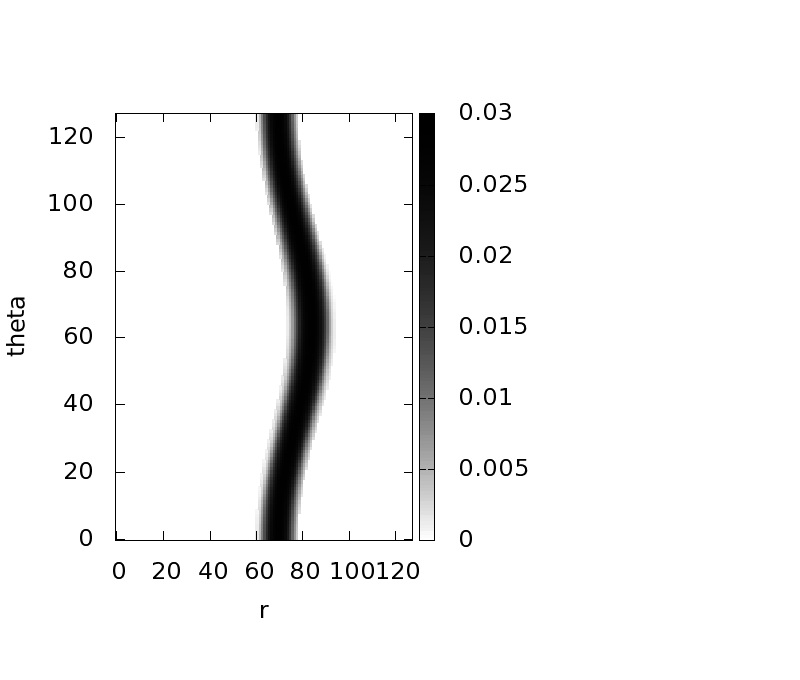}
\end{minipage}
\caption{Initial function $f_{t=0}$ with $v_{\parallel}\!=\!-3\,v_{th0}, N_r=128, N_\theta=128$ (left), distribution function $f$ at \mbox{$t=2000$} with FSL scheme (middle), distribution function $f$ at \mbox{$t=50000$} with FSL scheme (right)}
\label{shear10}
\end{figure}

\begin{figure}[h]
\begin{minipage}[t]{.48\textwidth}
\includegraphics[trim=25mm 22mm 80mm 55mm,clip,width=\textwidth]{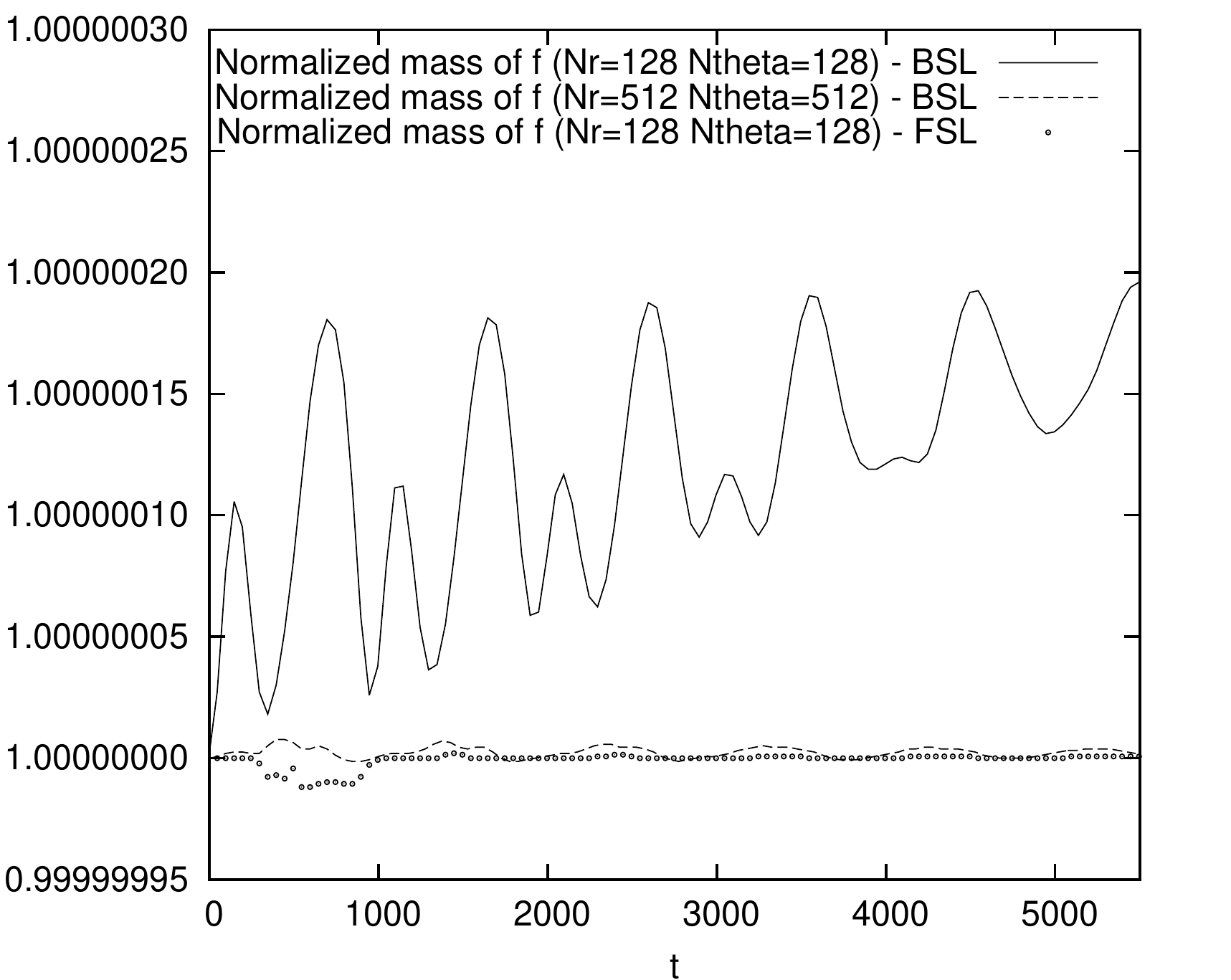}
\caption{Time evolution of the mass of $f$ for different poloidal mesh sizes}
\label{blob1}
\end{minipage}~
\begin{minipage}[t]{.48\textwidth}
\includegraphics[trim=25mm 22mm 80mm 55mm,clip,width=\textwidth]{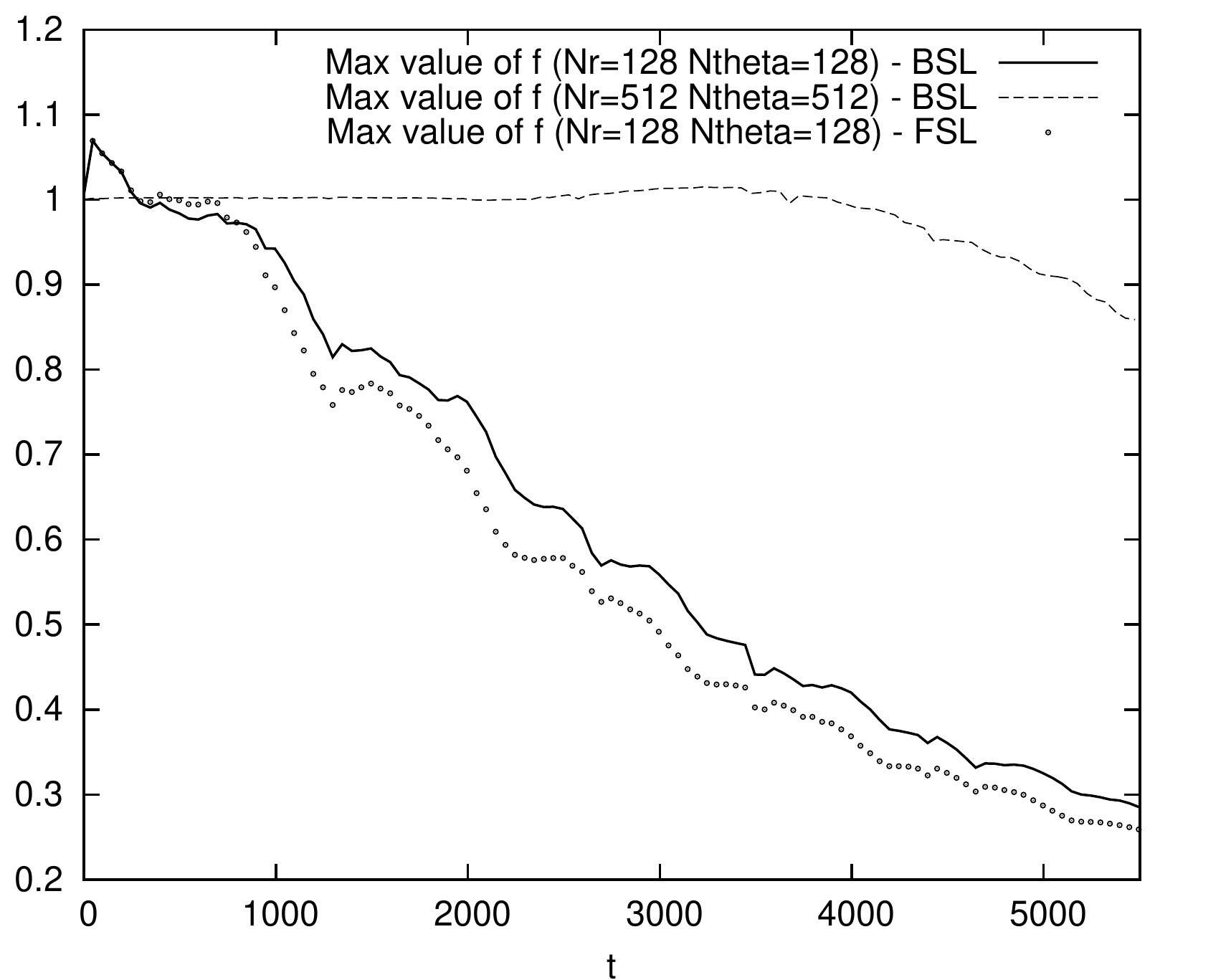}
\caption{Evolution of the $f$ maximum value for different poloidal mesh sizes}
\label{blob2}
\end{minipage}
\end{figure}

On Figures \ref{shear20}, \ref{shear0}, \ref{shear10}, the
distribution function is shown for different time steps (parallel
velocity is $v_{\parallel}\!=\!-3\,v_{th0}$, $m$ is set to 1, ,
$\rhostar=.01$). One can notice the shearing that appears along
time. The system dynamics tends to create filamentation in
configuration space. Nevertheless, the mass $\|f\|_1=\int f\,J_s\,
B_{\parallel}^{\star} dr d\theta d\phi dv_{\parallel}$ and the maximum
value $\|f\|_{\infty}=sup\{|f|:\forall
(r,\theta,\phi,v_{\parallel})\}$ should theoretically be conserved at
any time. Finer structures than mesh element size are dynamically
created (see the rightmost plot of Figure \ref{shear20}), and they can
not be correctly captured (filamentation induces that mass can be lost
or gained at each time step). As it can be seen in Figures \ref{blob1}
and \ref{blob2}, with the classical Backward Semi-Lagrangian scheme
$\|f\|_{\infty}$ quantity does not remain constant, $\|f\|_1$ is well
conserved at $(N_r=512,N_\theta=512)$ while $(N_r=128,N_\theta=128)$
does not lead to good conservation of $\|f\|_1$.  It is noteworthy,
that finer resolution in space - going from $(N_r=128,N_\theta=128)$
to $(N_r=512,N_\theta=512)$ greatly improves the conservation of mass
and of the largest absolute value. If one takes the Forward
Semi-Lagrangian scheme, the Figures \ref{blob1} and \ref{blob2}
illustrate that $\|f\|_{\infty}$ behaves like BSL but $\|f\|_1$ is
better preserved, as expected.
Then, FSL scheme is a good solution as far as mass
conservation is required.

\medskip
Let suppose one wants to keep the framework of the Backward
Semi-Lagrangian scheme and to avoid FSL, another solution to preserve
mass could be to have an aligned mesh. If one change the variables
$(r,\theta,\varphi,v_{\parallel})$ to
$(P_{\varphi},\theta,\varphi,v_{\parallel})$, the two initial
distribution functions we have seen are more accurately
discretized. Also, the conservation of mass is improved for the
following reason: advection is done only in the direction $\theta$,
there is no shift along the $P_{\varphi}$ direction. The displacement
in the poloidal plane during one time step amounts to shift only in
$\theta$. So mass has just to be conserved on each mesh line labelled
by $P_{\varphi}$, and it will automatically give global mass
conservation. It is easy to construct an advection operator and the
associated interpolation scheme that preserve mass on each mesh line
$P_{\varphi}$. Another benefit is that in this new set of variables,
the equilibrium function is better discretized and represented: steep
gradient is only in $P_{\varphi}$ direction in the poloidal
plane. Nevertheless, building and managing
a mesh that includes $P_{\varphi}$ variable is much more difficult to
implement.


\subsubsection{Extrapolation to nonlinear test cases}

The first test case of section \ref{firstcase} is linked to what is
going on at the beginning of a Gysela run. The initial distribution
function is a perturbed equilibrium. One wishes to remain close to
this quasi-equilibrium (if the amplitude of electric potential is very
small) without beeing too much biased by some large numerical
artifacts of the linear operator that we have seen: artificial
increase of $u_1$, $u_\infty$.

The second test case of section \ref{secondcase} is near what is
happening in non-linear dynamics in which many localized spatial
structures develop and grow. One expect to conserve the total particle
density, even if the distribution function becomes something else than
an ideal function of motion invariants. The decrease of maximal value
of Figure \ref{blob2} is not acceptable to perform accurate
simulations. Aligned coordinates along $P_{\varphi}$ should be able to
improve conservation of $\|f\|_{\infty}$.

\section{Nonlinear gyrokinetic simulations}

\subsection{Split linear and nonlinear parts}

\subsubsection{Global separation of linear/nonlinear terms}
The equations (\ref{dxidt},\ref{dvpardt}) can be split into two parts
(in the same spirit that \cite{Ido08}). The first part includes the
nonlinear terms that depends on the electric potential. The second
part comprises the other terms. One can solve these two parts one
after the other by splitting. On the first hand, the linear operator
is presented in Eqs~(\ref{lin1},\ref{lin2}), on the second hand, the
nonlinear operator is described by Eqs~(\ref{nl1},\ref{nl2}).\\

\begin{minipage}{5.9cm}
\centering
\textbf{Linear operator $\mathcal{L}$}
\footnotesize
\begin{align}
\frac{{\rm d}x^i}{\rm dt} &= \vpar \vecbstar\cdot\gradxi + \vDgradxi \label{lin1}\\
m_s\frac{{\rm d}\vpar}{\rm dt} &= -\mu \vecbstar\cdot\gradB \label{lin2}
\end{align}
\end{minipage}~
\begin{minipage}{7.6cm}
\centering
\textbf{Nonlinear operator $\mathcal{N}$}
\footnotesize
\begin{align}
  \frac{{\rm d}x^i}{\rm dt} &= \vgcEgradxi \label{nl1}\\
  m_s\frac{{\rm d}\vpar}{\rm dt} &= - e_s \vecbstar\cdot\gradgyroPhi + \frac{m_s\vpar}{B}\vgcEgradB \label{nl2}
\end{align}
\end{minipage}

\vspace*{.7cm} The linear operator exhibits large displacements at
large modulus of parallel velocity, and also induces shear
flows. These features can interfere with the nonlinear dynamics that
eventually involves small displacements that are not of the same order
of magnitude. Moreover, as the dynamics generated by the two operators
are different, induced accuracy problems have possibly not the same
characteristics for the two operators; and the limitations (CFL-like
conditions) on the time step are also not the same. Ideally, one
should be able to fix the time step of linear and nonlinear operators
independantly in order to enhance accuracy.

The original formulation that uses Strang splitting in phase space for
the Vlasov solver, but without a linear/nonlinear splitting, can leed
to serious troubles at high parallel velocities. Indeed, directional
splitting suffers from some large shifts in $\varphi, r$ and $\theta$
at each directional advection (section \ref{vlasovsplit}). Then, for
not so large $\Delta{}t$, the evaluation of electric fields $E_{(.)}$
is not done at the right spatial position at each substep of the
splitting except for the first substep of the splitting. One has to
take very small time step to recover small displacements in the linear
operator and then a reasonable accuracy in the evaluation of
$E_{(.)}$.\\ The splitting between linear and nonlinear parts
eliminates this problem. The nonlinear operator is applied alone, thus
the linear operator and its large shifts does not interact badly with
nonlinear solver.  The proposed solution is to perform at each time
step: first, a directional splitting for the linear operator
$\mathcal{L}$ with the sequence
$(\hat{\vpar}/2,\hat{\varphi}/2,\hat{r\theta},\hat{\varphi}/2,\hat{\vpar}/2)$,
second, a directional splitting for the nonlinear operator
$\mathcal{N}$ with the same sequence
$(\hat{\vpar}/2,\hat{\varphi}/2,\hat{r\theta},\hat{\varphi}/2,\hat{\vpar}/2)$.


This approach is a liitle bit expensive in term of computational
cost because it nearly doubles time to solution. 


\subsubsection{Conservation of mass}

Each of the operators $\mathcal{L}$, $\mathcal{N}$ conserves the mass
independently. One can track during simulation which operator degrades
or not the total mass.

\subsection{Experiments}

In the following GYSELA runs, the BSL and FSL schemes are compared for
two test cases with $\mu$ =0 (4D only).
Results for both BSL scheme and FSL scheme
are shown. In
the FSL version that has been run, FSL scheme is present in all
advections: 1D ($\varphi,\vpar$) and 2D ($r,\theta$). The rest of the code is unchanged compared
to the BSL version. In order to simplify boundary conditions in
$\vpar$, the distribution function is fixed at $v_{min}$ and
$v_{max}$. Therefore, we have also taken this choice for the BSL runs presented hereafter.

\subsubsection{Cylindrical case}

\begin{figure}[H]
\begin{minipage}[t]{.45\textwidth}
\includegraphics[trim=19mm 22mm 8mm 20mm,clip,width=\textwidth]{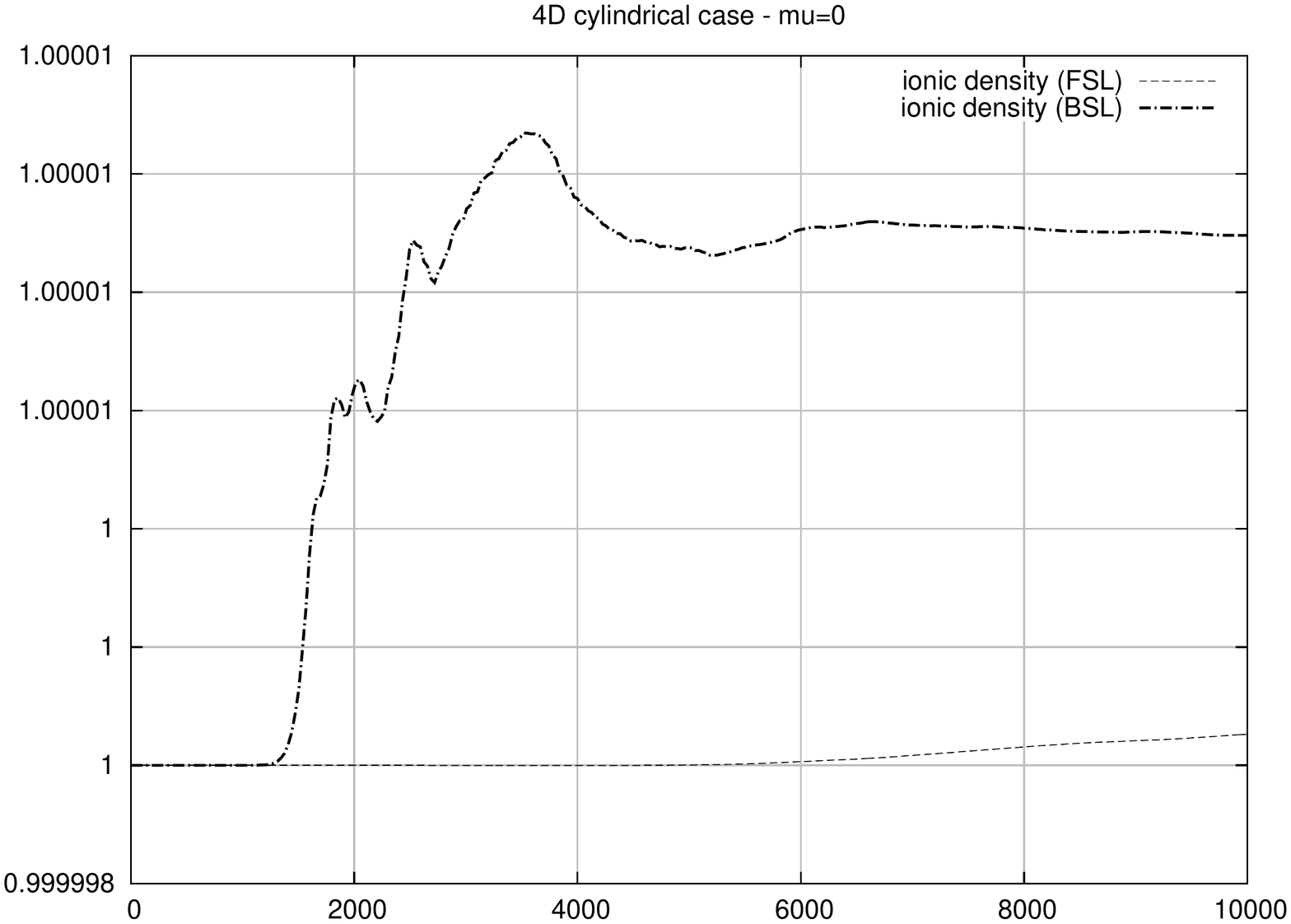}
\caption{Evolution of the mass for a 4D, $\mu=0$, cylindrical test case}
\label{cyl1}
\end{minipage}~
\begin{minipage}[t]{.45\textwidth}
\includegraphics[trim=19mm 22mm 8mm 20mm,clip,width=\textwidth]{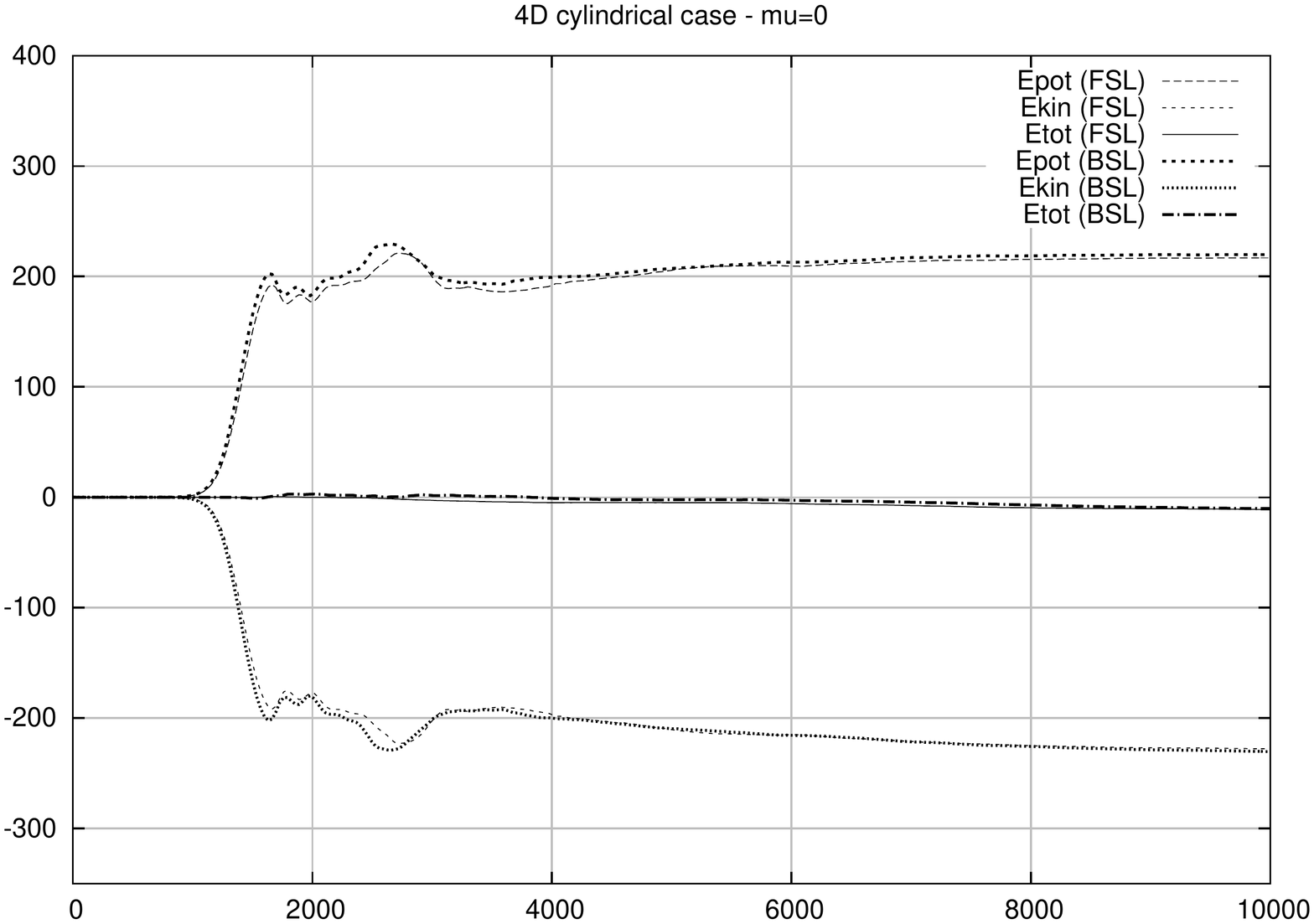}
\caption{Evolution of the energies for a 4D, $\mu=0$, cylindrical test case}
\label{cyl2}
\end{minipage}
\end{figure}

Some plots for our reference test case \texttt{data4D\_slab\_A32}
($\rhostar=1/32$, $\mu=0$) is given in
Fig.~(\ref{cyl1},\ref{cyl2},\ref{cyl3},\ref{cyl4}).  The computational
domain size is $N_r=128,N_\theta=256,N_\varphi=32,N_{\vpar}=48$.

The mass is well preserved with FSL as it can be seen on Fig.~\ref{cyl1}. Potential energy, kinetic energy and total energy are presented on Fig.\ref{cyl2}. The total energy should theoretically remains at zero, the two schemes are quite proximate to that. During the beginning of the simulation, FSL and BSL give quite similar electric potential as it is shown at time $t=1920$ in Fig.~(\ref{cyl3},\ref{cyl4}). Nevertheless, after a while, the two schemes diverges slowly (see for example $t=8000$).

\begin{figure}[H]
\begin{minipage}[t]{.49\textwidth}
\includegraphics[trim=10mm 15mm 4mm 5mm,clip,width=\textwidth]{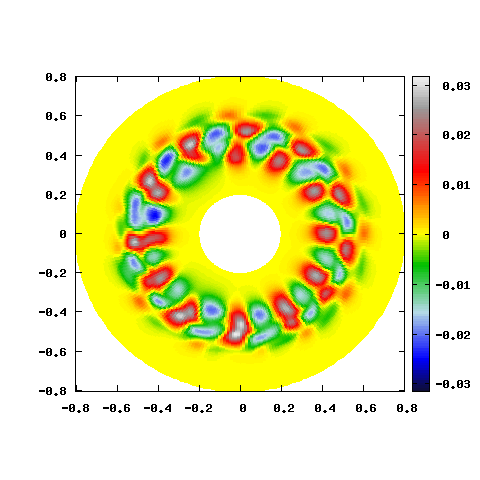}
\end{minipage}
\begin{minipage}[t]{.49\textwidth}
\includegraphics[trim=10mm 15mm 4mm 5mm,clip,width=\textwidth]{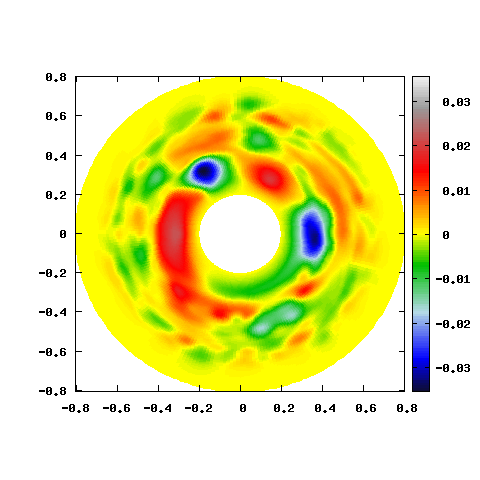} 
\end{minipage}
\caption{Poloidal cut of electric potential at different time steps ($t=1920$ left, $t=8000$ right), $\varphi=0$, \textbf{BSL} scheme}
\label{cyl3}
\begin{minipage}[t]{.49\textwidth}
\includegraphics[trim=10mm 15mm 4mm 5mm,clip,width=\textwidth]{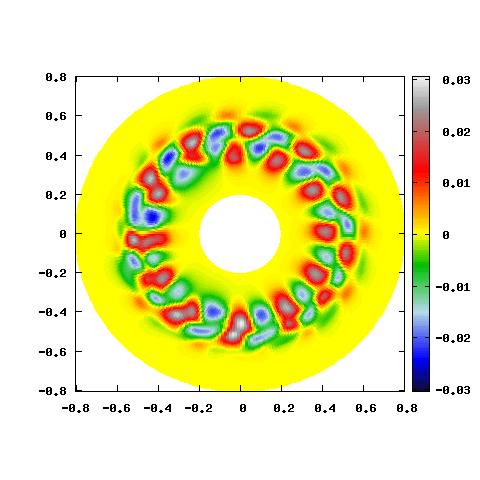}
\end{minipage}
\begin{minipage}[t]{.49\textwidth}
\includegraphics[trim=10mm 15mm 4mm 5mm,clip,width=\textwidth]{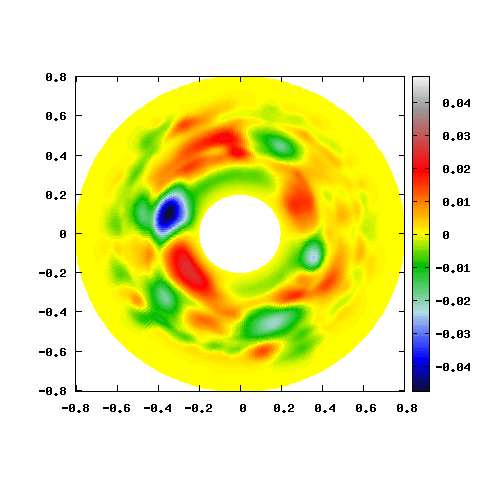} 
\end{minipage}
\caption{Poloidal cut of electric potential at different time steps ($t=1920$ left, $t=8000$ right), $\varphi=0$, \textbf{FSL} scheme}
\label{cyl4}
\end{figure}

\subsubsection{Toroidal case}

A small toroidal test case has been set up ($\rhostar=1/40$).  The
computational domain size is
$N_r=256,N_\theta=256,N_\varphi=64,N_{\vpar}=48$.  The FSL version
outperforms the BSL version from the point of view of mass
conservation. The conservation of total energy is far more better with
FSL and this is due to improved mass conservation. A precise tracking
of mass losses/gains allows us to conclude that the linear operator
$\mathcal{L}$ is mainly responsible for mass degradation using BSL
approach. An explanation is that in toroidal setting, $\mathcal{L}$
causes large displacements, whereas in cylindrical setting (previous
subsection) this was not the case.

\begin{figure}[H]
\begin{minipage}[t]{.98\textwidth}
\centering
\includegraphics[trim=9mm 20mm 7mm 15mm,clip,width=.5\textwidth]{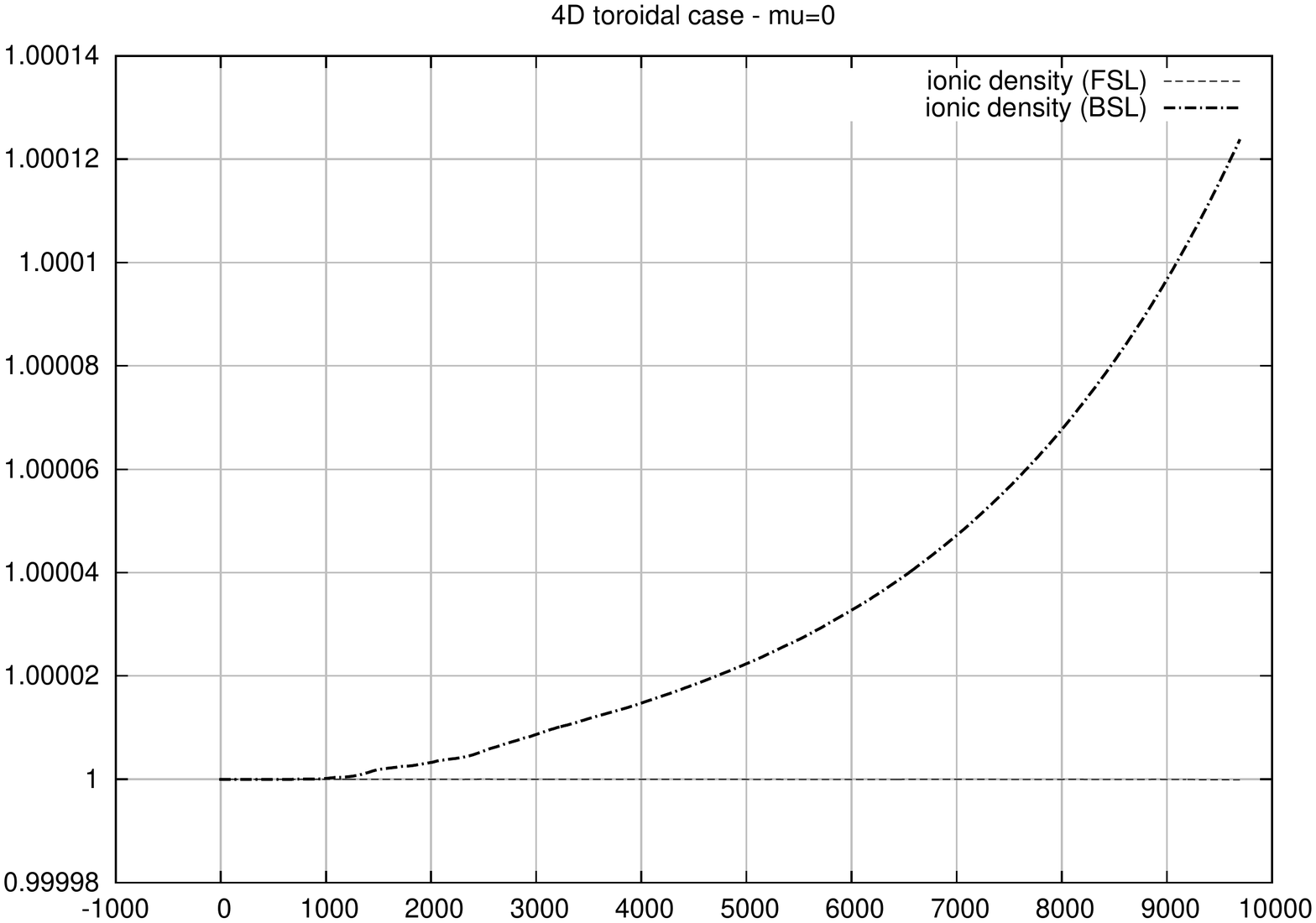}
\caption{Evolution of the mass for a 4D, $\mu=0$, toroidal test case}
\label{tor1}
\end{minipage}
\begin{minipage}[t]{.45\textwidth}
\includegraphics[trim=9mm 20mm 7mm 2mm,clip,width=\textwidth]{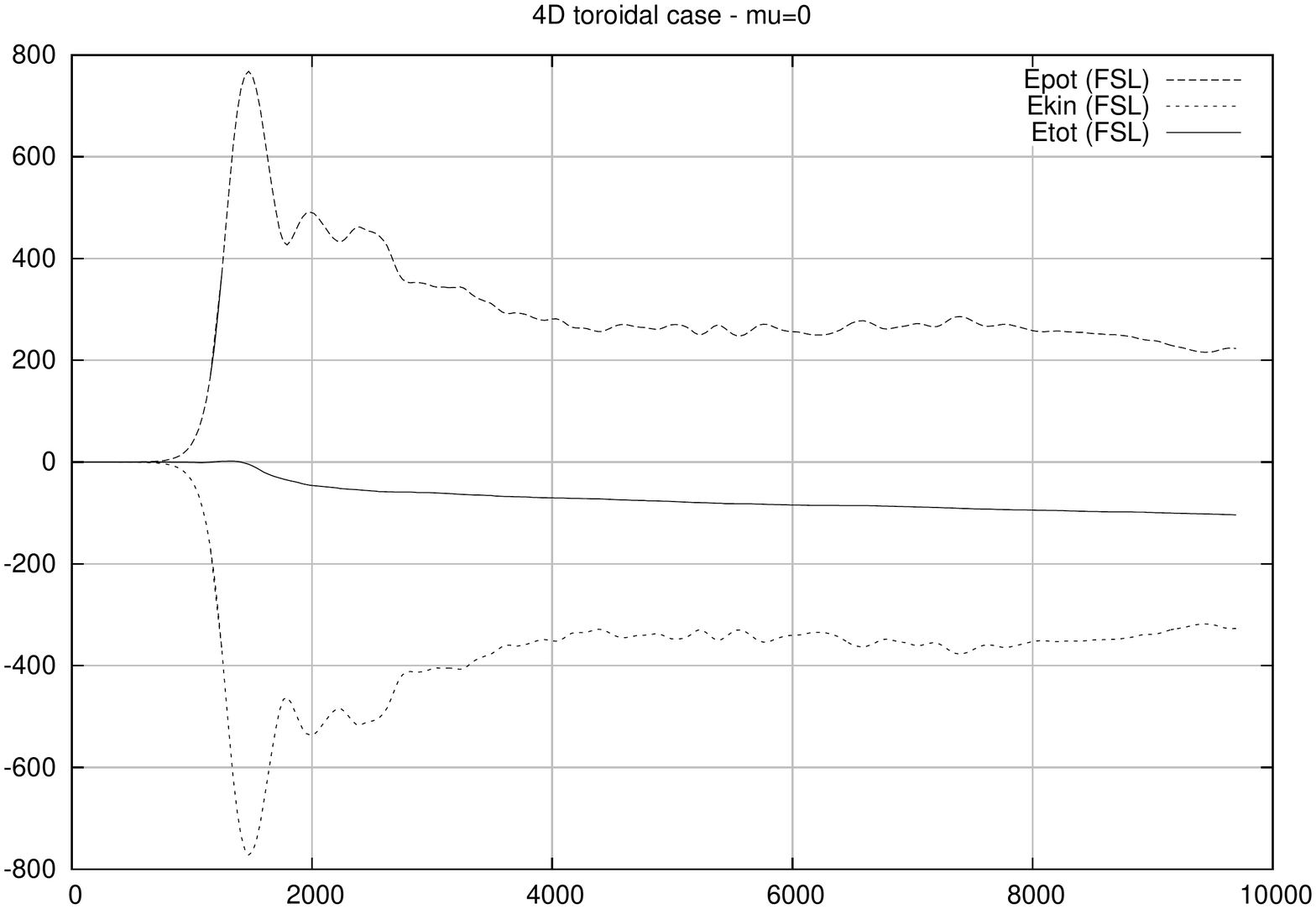}
\caption{Evolution of the energies for a 4D, $\mu=0$, toroidal test case with FSL}
\label{tor2}
\end{minipage}~
\begin{minipage}[t]{.45\textwidth}
\includegraphics[trim=9mm 20mm 7mm 2mm,clip,width=\textwidth]{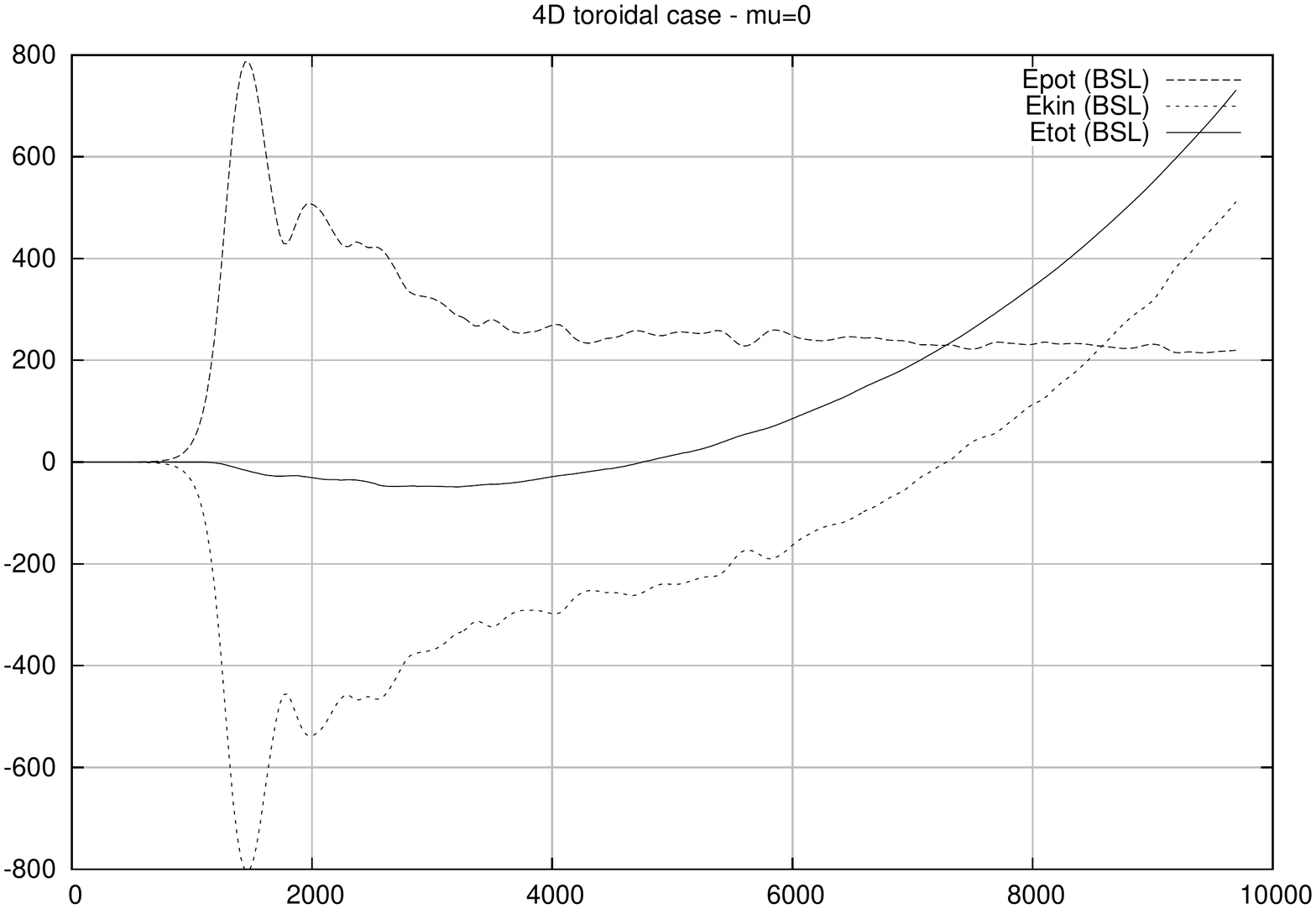}
\caption{Evolution of the energies for a 4D, $\mu=0$, toroidal test case with BSL}
\label{tor3}
\end{minipage}
\begin{minipage}[t]{.42\textwidth}
\includegraphics[trim=10mm 15mm 4mm 5mm,clip,width=\textwidth]{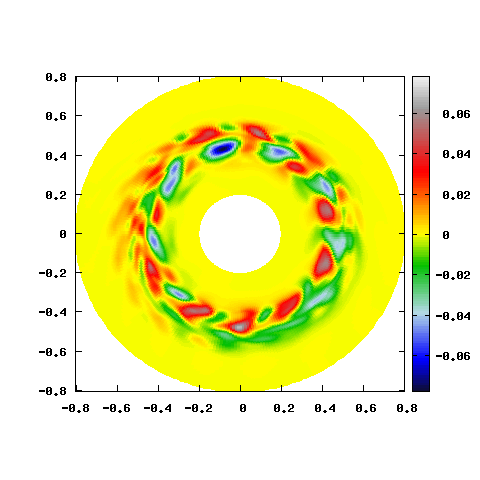}
\end{minipage}\hspace*{1.8cm}
\begin{minipage}[t]{.42\textwidth}
\includegraphics[trim=10mm 15mm 4mm 5mm,clip,width=\textwidth]{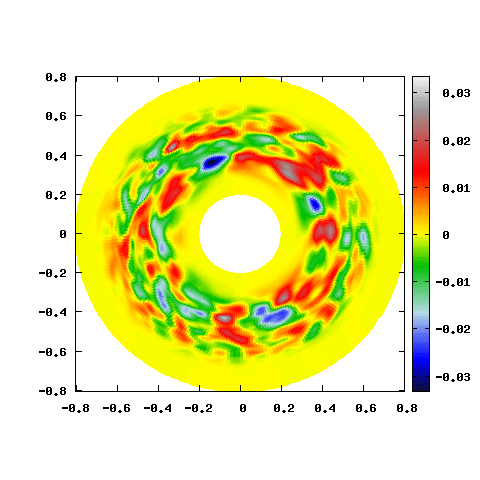} 
\end{minipage}
\caption{Poloidal cut of electric potential at different time steps ($t=1600$ left, $t=4800$ right), $\varphi=0$, \textbf{BSL} scheme}
\label{cyl3}
\begin{minipage}[t]{.42\textwidth}
\includegraphics[trim=10mm 15mm 4mm 5mm,clip,width=\textwidth]{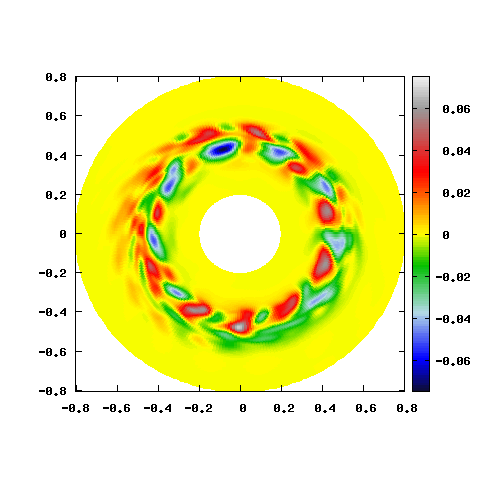}
\end{minipage}\hspace*{1.8cm}
\begin{minipage}[t]{.42\textwidth}
\includegraphics[trim=10mm 15mm 4mm 5mm,clip,width=\textwidth]{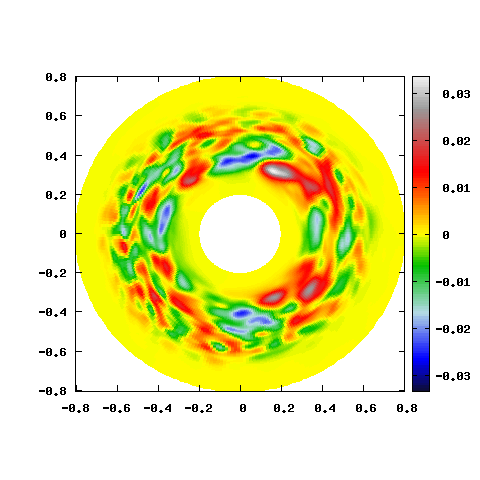} 
\end{minipage}
\caption{Poloidal cut of electric potential at different time steps ($t=1600$ left, $t=4800$ right), $\varphi=0$, \textbf{FSL} scheme}
\label{cyl4}
\end{figure}




\section{Conclusion}

Few solutions have been proposed to help preserving constant states,
$L_1$ and $\L_{\infty}$ norms in gyrokinetic semi-Lagrangian
simulations. Some GYSELA runs show that we get benefits from these
improvements in term of accuracy.

\bibliographystyle{alpha}
\bibliography{eqnInvar}

\end{document}